\documentclass{article}

\usepackage{amsmath}
\usepackage{graphicx}
\newtheorem{theorem}{Theorem}
\newtheorem{lemma}{Lemma}
\newtheorem{proposition}{Proposition}
\newtheorem{corollary}{Corollary}
\newtheorem{definition}{Definition}

\title{Linear Statistics of Non-Hermitian Matrices Matching the Real or Complex Ginibre 
\\ Ensemble to Four Moments}
\author{Phil Kopel}
\date{}

\begin{document}
\maketitle

\begin{abstract}
We prove that, for general test functions, the limiting behavior of the linear statistic of an independent entry random matrix is determined
only by the first four moments of the entry distributions. This immediately generalizes the known central limit theorem for independent entry matrices with complex normal entries.
We also establish two central limit theorems for matrices with real normal entries, considering seperately functions supported exclusively on and exclusively away from the real line. In contrast to previously obtained results in this area, we do not impose analyticity on test functions.
\end{abstract}

\section{Introduction}

This paper deals with linear statistics of random matrices. These are the random matrix analogues of the central limit theorems of classical probability theory,
and describe the fluctuations of the spectra of random matrices about their expected configurations. As is common in random matrix theory, we will be concerned with the limiting behavior of these statistics as the matrix dimension becomes large.

The general framework is as follows: let $f$ be some test function and $M$  some $n\times n$ random matrix with eigenvalues $\lambda_j$. One then defines the associated \textit{linear statistic} to be:

\begin{eqnarray}
S_n[f]=\sum_{j=1}^n f\left(\lambda_j\right)-\mathbf{E}\left[
\sum_{j=1}^n f\left(\lambda_j\right)
\right]
\end{eqnarray}

Depending on the random matrix ensemble and test function under consideration, this quantity is sometimes normalized and sometimes not -- we will encounter both situations. This statistic is understood for many models, most famously the well known Wigner random matrix ensembles (Hermitian matrices with independent entries up to the Hermiticity condition), whose unnormalized linear statistic is known to converge to a normal distribution for a wide range of test functions \cite{SW}.

The random matrices that we will consider here will be \textit{independent entry random matrices}, which are sometimes simply referred to simply as non-Hermitian random matrices. These are the non-Hermitian analogues of
the classical Wigner random matrices, and their construction is straightforward: each entry is an independent random variable. These matrices find application in physics, neuroscience, economics, and other disciplines (the interested reader is directed to \cite{AK} and the references therein). The precise definition we will use is:

\begin{definition}
Let $\xi_{i,j}$ be independent random variables with mean zero and variance one,
which match one another to four moments. Assume further that for all $t >0$:
\begin{eqnarray*}
\mathbf{P}(|\xi_{i,j}|>t)\leq Ce^{-t^c}
\end{eqnarray*}

Then the random matrix $M_n=(n^{-1/2}\xi_{i,j})_{1\leq i,j\leq n}$ is said to be an \textit{independent entry random matrix}. The individual entries $\xi_{i,j}$ are called the \textit{atom distributions}.
\end{definition}

Less is known about independent entry matrices than about their Hermitian counterparts, and their study is largely more recent. There are a number of good reasons for this -- for one thing,
polynomials are no longer dense among test functions, which limits the effectiveness of several methods.  Eigenvalues stability inequalities also fail in the non-Hermitian case:
two non-Hermitian matrices may look almost identical but have radically different spectral properities. 

Nevertheless, we have a fairly good understanding of what the eigenvalues of these ensembles look like. The celebrated \textit{circular law} states that the limiting distribution
of the eigenvalues of an independent entry matrix is uniform on the unit disc $D(0,1)$. Informally, for nice functions $f$ and large independent entry matrices, the circular law states:

\begin{eqnarray}
\frac{1}{n}\sum_{j=1}^n f(\lambda_j)\approx \frac{1}{2\pi}\int_{D(0,1)} f(z)d^2z
\end{eqnarray}

This analogue of the Wigner semi-circular law was originally proven in the case of complex Gaussian entries by Ginibre in 1965, this ensemble now is known as the complex \textit{Ginibre ensemble}. (The complex Ginibre ensemble is sometimes also referred to as the GinUE, and the real Ginibre ensemble as the GinOE, in analogy to the Wigner GUE and GOE ensembles.) It was not immediately obvious how to extend Ginibre's result to more general independent entry matrices, and indeed the circular law was only proven in full generality in recent years \cite{TV1} -- see the survey \cite{BD}.

The study of fluctuations about the circular law appears to have been initiated in a 1999 paper of Forrester \cite{For},
which considers the case of complex normal entries and radially invariant test functions. The calculation essentially reduces to the orthogonality of monic polynomials with
respect to radially invariant measures, and a normal limiting distribution is obtained for the unnormalized linear statistic.

Forrester's result has been extended considerably by a pair of papers, the first by Rider and Virag \cite{RV} and the second by Rider and Silverstein \cite{RS},
which establish central limit theorems in the case of complex normal entries and general test functions, and in the case of general entries but analytic test functions, respectively. The first employs a combinatorial proof which relies on explicit formulas available in the case of complex Gaussian entries,
while the second uses a martingale difference argument.

It remains, then, to investigate situations where test functions need not be analytic and matrix entries need not have complex normal distributions. The main result presented here is a universality principle: we show that
the limiting distribution of the linear statistic only depends on the first four moments of the atom distributions. This behavior is exactly analagous to that of the linear statistics of Wigner matrices, whose limits are also only sensitive to the first four moments of the atoms. 

Colloquially, this means that if the linear statistic is understood for one matrix model, it is also understood for sufficiently similar
models. The precise statement is:

\begin{theorem}
\label{nh:4mom}
Let $M_n^1$  and $M_n^2$ be two independent entry random matrices whose atom distributions match to four moments, and let $f$ be a smooth function with compact support. If the linear statistic for
$M_n^1$, denoted $S^1_n[f]$, converges weakly to some limiting distribution $\chi$ as $n \to \infty$, then  the linear statistic associated with $M_n^2$, denoted $S_n^2[f]$, converges weakly to $\chi$ as well.
\end{theorem}

In fact, we will prove a somewhat more general four moment lemma, which will imply Theorem 1. The proof will use the same resolvent swapping technique that was successfully applied in \cite{TV2}, in combination with the recent eigenvalue rigidity estimates obtained in \cite{BYY}.

Comparing with Rider and Virag's result, we immediately obtain the following corollary:

\begin{corollary}
\label{nh:rvmatch}
Let $M_n$ be an i.i.d. random matrix whose entries have mean zero and unit variance. Suppose the first four moments of these distibutions match the complex Gaussian.
 Let $f$ be a smooth test function, and let $D(0,1)$ denote the unit disc.

Then $N_n[f]$ converges weakly to a normal random variable with variance $\sigma^2_A
+\sigma^2_B$, where:
\begin{eqnarray*}
\sigma^2_A=\frac{1}{4}||\phi||^2_{H^1(\mathbf{D(0,1)})}
\end{eqnarray*}

And:
\begin{eqnarray*}
\sigma^2_B=\frac{1}{2}||f||^2_{H^{1/2}(\mathbf{\delta D(0,1)})}
\end{eqnarray*}

These are known as the bulk component and the edge component, respectively, of  the variance.
\end{corollary}

In addition to the base case of complex Gaussian entries, we will also consider real Gaussian entries. Real entries force the spectrum to split into a real and complex components: the density of the eigenvalues
is not dominated by the two-dimensional Lebesgue measure on the real line. For this reason, we will consider functions supported either
entirely on or entirely away from the real line seperately (and since eigenvalues come in conjugate pairs, we will work exclusively in the upper half plane without loss of generality).

Away from the real line, we will obtain:

\begin{theorem}
\label{nh:bulk}
Let $M_{2n}$ be an independent entry random matrix whose entries match real Gaussians to four moments. If $f$ is a smooth test function compactly supported in the upper half of the unit disc away from the real line,
then the linear statistic $S_{2n}[f]$ converges in distribution to a normal random variable with variance:
\begin{eqnarray*}
\frac{1}{4\pi}\int \left(f_x(x,y)^2+f_y(x,y)^2\right) dxdy
\end{eqnarray*}
\end{theorem}

On the real line, we will obtain:

\begin{theorem}
\label{nh:line}
Let $M_{2n}$ be an independent entry random matrix whose entries match real Gaussians to four moments. If $f$ is a smooth test function compactly supported on the interval $(-1,1)$,
then the normalized linear statistic $n^{-1/4}S_{2n}[f]$ converges in distribution to a normal random variable with variance:
\begin{eqnarray*}
\left(\frac{2-\sqrt{2}}{\sqrt{\pi}}\right)\int f(x)^2 dx
\end{eqnarray*}
\end{theorem}

These results will follow from explicit formulas available for these matrices, as well as some useful quaternion identities, which we will present in due course. The case of test functions supported away from the real line is substantially more computationally intensive to prove, because it is necessary to show that all the higher cumulants will cancel out in the limit. In the case of test functions supported on the real line, on the other hand, we merely have to show that the normalization suppresses all higher cumulants -- no cancellation is required.

There are several potential directions for future work. First, one may consider the real Ginibre ensemble and test functions whose support intersects but is not contained in the real line. Second, one may consider the real Ginibre ensemble and functions whose support extends to the edge of the spectral support (so to the circle $|z|^2=1$). The methods applied here should basically generalize to these cases, although we do not pursue these directions here. Most importantly, it remains to find central limit theorems for independent entry matrices which do not match the real or complex Ginibre ensemble to four moments, without resorting to the assumption of analytic test functions. As the methods employed here for the Ginibre ensemble depend very crucially on explicit formulas not available for matrices with non-Gaussian entries, our methods are unlikely to generalize to this problem, and new ideas are needed.

The paper is organized as follows: in Section 2, we will first introduce relevant techniques and results, and then prove Theorem \ref{nh:4mom} (which as noted previously, immediately implies Corollary 1). In Section 3, we will again first introduce relevant result (and some of the basics of the theory of quaternions), then we will prove Theorem \ref{nh:bulk} and Theorem \ref{nh:line} in that order.

\section{Universality and Central Limit Theorem for GinUE}

\subsection{Preliminary Reductions}

In this section, we introduce a few techniques and known results in preperation for the proof of Theorem \ref{nh:4mom}. We begin by rewriting the linear statistic $S_n[f]$ in a more tractable form.
 For any complex number $z$ and matrix $M_n$, define the Hermitian matrix $W_n(z)$ by:

\begin{eqnarray}
\label{girkoh}
W_n(z)=\frac{1}{\sqrt{n}}\left(
\begin{array} {cc}
0 & (M_n-z)\\
(M_n-z)^{*} & 0
\end{array}
\right)
\end{eqnarray}

We will occasionally take the Stiejles transform of this matrix:

\begin{eqnarray}
s(z,\zeta)=\frac{1}{n}\mbox{Tr}(W_{n}(z)-\zeta)^{-1}
\end{eqnarray}

Suppose $f: \mathcal{R}^2 \to \mathcal{R}$ is a twice differentiable function, compactly supported on the ball $B(0,1-\epsilon)$, for some small $\epsilon > 0$. Then by the Cauchy-Pompeiu formula and integration by parts:
\begin{eqnarray*}
\sum_{j=1}^n f(\lambda_j)=\sum_{j=1}^n\frac{-1}{2\pi}\int_{B(0,1)}\frac{df(z)}{d\bar{z}}
\frac{1}{z-\lambda_j}dxdy\\
=\frac{1}{2\pi}\int_{B(0,1)}\Delta f(x,y)\log(|\prod_j((x+iy)-\lambda_j)|)dxdy\\
+\frac{i}{2\pi}\int_{B(0,1)}\Delta f(x,y)\mbox{atan2}
(\Re(\prod_j((x+iy)-\lambda_j)),\Im(\prod_j((x+iy)-\lambda_j)))dxdy
\end{eqnarray*}

Since the first expression is real then the second must be as well, and we conclude:
\begin{eqnarray*}
\sum_{j=1}^n f(\lambda_j)=\frac{1}{2\pi}
\int_{B(0,1)}\Delta f(z)\log(|\det(W_{n,z})|)d^2z
\end{eqnarray*}

Or else:
\begin{eqnarray}
S_n[f]=\frac{1}{8\pi}\int_{|z|\leq 1}\Delta f(z)\log|\det(M_n-z)|d^2z
\end{eqnarray}

This reformulation is commonly called Girko's Hermitization trick, and is an essential component of the proof of the general circular law \cite{BD}.

Next, we present three lemmas which feature in the proof.  The following Monte Carlo Sampling lemma follows from an easy application of the Chebycheff inequality \cite{TV2}:

\begin{lemma}\label{p1}
Let $(X,\mu)$ be a probability space and $F$ a square integrable function from $(X,\mu)$ to the real line. For $m$ independent, uniformly distributed $x_i$  define the empirical average:
\begin{eqnarray*}
S_m=\frac{1}{m}\sum_{i=1}^m F(x_i)
\end{eqnarray*}
Then for any $\delta >0$ the following estimate holds with probability at least $1-\delta$:
\begin{eqnarray}
|S_m-\int_X F d\mu|\leq \frac{1}{\sqrt{\delta m}}\left(\int_X (F-\int_X Fd\mu)^2d\mu\right)^{1/2}
\end{eqnarray}
This is a probabilistic variant of estimating an integral by a Riemann sum
\end{lemma}

Notice that we may rewrite this estimate in the following, occasionally more useful, form:
\begin{eqnarray}
\mathbf{P}\left[ |S_m-\int_X F d\mu|\geq \delta \right]\leq \frac{\int_X (F-\int_X Fd\mu)^2d\mu}{\delta^2 M}
\end{eqnarray}

In addition to Monte Carlo sampling, we will also employ the technique of resolvent swapping, which is a probabilistic analogue of Taylor expansion.

Define an \textit{elementary matrix} to be a Hermitian matrix featuring only one or two entries equal to a unit multiple of 1, and all the other entries set to zero.  Therefore, adding a multiple of an elementary matrix to a Hermitian matrix $M$
changes either a single diagonal entry or two conjugate off-diagonal entries of $M$, and leaves the other entries alone.

For a Hermitian matrix $M$ and an elementary matrix $V$, define:
\begin{eqnarray*}
M_t=M+\frac{1}{\sqrt{n}}tV\\
R_t(\zeta)=(M_t-\zeta)^{-1}\\
s_t(\zeta)=\frac{1}{n}\mbox{Tr}R_t(\zeta)
\end{eqnarray*}

We will also need to define an appropriate matrix norm:
\begin{eqnarray*}
||A||_{\infty,1}=\max_{1\leq i,j\leq n}|A_{ij}|
\end{eqnarray*}

The following Taylor expansion type lemma is due to Tao and Vu \cite{TV3}, and is proven by iterating the classical resolvent identity:
\begin{lemma}\label{p2}
Let $M_0$ be a Hermitian matrix, $V$ an elementary matrix, $t$ a real number. Let $R_0=M_0^{-1}$. Suppose we have:
\begin{eqnarray*}
|t|\times||R_0||_{(\infty,1)}=o(\sqrt{n})
\end{eqnarray*}
Then we have the following Taylor expansion:
\begin{eqnarray*}
s_t=s_0+\sum_{j=1}^k n^{-j/2}c_jt^j+O(n^{-(k+1)/2}|t|^{k+1}||R_0||_{(\infty,1)}
\min(||R_0||_{(\infty,1)},\frac{1}{n\eta})
\end{eqnarray*}
The coefficients $c_j$ are independent of $t$ and obey the following estimate:
\begin{eqnarray*}|c_j|\leq
||R_0||_{(\infty,1)}^j\min(||R_0||_{(\infty,1)},\frac{1}{n\eta})
\end{eqnarray*}
\end{lemma}

For a proof, see \cite{TV3}.

Finally, we will need a rigidity estimate for the eigenvalues of an independent entry matrix. Define $p_c(w,z)$ to be the function such that $\int_{\mathbf{R}}\frac{p_c(x,z)}{x-w}dz=m_c(w,z)$, where $m_c(w,z)$ is the unique solution to the following equation:
\begin{eqnarray*} 
m_c^{-1}=-w(1+m_c)+|z|^2(1+m_c)^{-1}
\end{eqnarray*}

Then classical position of the $j$-th eigenvalue, $\gamma_j(z)$, is given by the formula:
\begin{eqnarray*}
\int_0^{\gamma_{j}(z)}p_c(x,z)dx=j/N
\end{eqnarray*}

It is unlikely for the eigenvalues of a random independent entry matrix to be located very far from their classical positions. Yau, Yin and Bourgade \cite{BYY} have shown:

\begin{lemma}\label{p3}
If $\lambda_j$ are the eigenvalues of an independent entry matrix, the following estimate holds uniformly for $z$ a fixed distance away from the spectral edge $|z|=1$:
\begin{eqnarray*}
\mathbf{P}\left[\left|\sum_{j=1}^n \left(\log\lambda_j(z)-\log\gamma_j(z)\right)\right|^2
\geq \log(n)^{C\log\log(n)}\right]\leq n^Ce^{-\log(n)^{C\log\log(n)}}
\end{eqnarray*}
In particular, this holds with overwhelming probability.
\end{lemma}

This is a consequence of the proof of Lemma 2.2 in \cite{BYY}.

\subsection{Proof of Theorem \ref{nh:4mom}}

Here, we prove the following four moment lemma:

\begin{lemma}\label{4m}
Let $G(t)$ be a bounded function with at least five bounded derivatives. Let $M_n^{1}$ and $M_n^2$ be two independent entry matrices, with atom distributions matching to four moments. Let $f$ be a twice continuously differentiable function compactly supported in the interior of the unit disc. Then:
\begin{eqnarray*}
\mathbf{E}\Big[G\left(\int \Delta f(z)\log|\det(W^1_n(z))|dz\right)-
G\left(\int \Delta f(z)\log|\det(W^2_n(z))|dz\right)\Big]
\\
=O(n^{-1/2+c_0})
\end{eqnarray*}
Additonally:
\begin{eqnarray*}
G\left(\mathbf{E}\int \Delta f(z)\log|\det(W^1_n(z))|dz\right)-
G\left(\mathbf{E}\int \Delta f(z)\log|\det(W^2_n(z))|dz\right)
\\
=O(n^{-1/2+c_0})
\end{eqnarray*}
Here, $c_0$ is a sufficiently small positive constant.
\end{lemma}

The proof is an adaptation of the proof of the similar Four Moments Lemma in \cite{TV2}. Notice that applying this lemma with $G=e^{-ixt}$ (for arbitrary $x$) proves Theorem \ref{nh:4mom} by way of the Fourier inversion theorem.

\textbf{Proof:}

We focus now on proving the first estimate in lemma \ref{4m}, with the second proven similarly. Notice that since we are studying the expectation of a bounded function we can (and will) treat any event which occurs with probability approaching one as occurring surely.

Let $M_n^1$ and $M_n^2$ be two independent entry matrices whose entries match to four moments, as in the statement of the theorem, and define the matrices $W^k_n(z)$ using Girko's Hermitization method, for $k=1,2$. Let $\lambda_j^k(z)$ denote the $k$-th eigenvalues of $W_n^k(z)$,

For $\beta=1,2$ and $m$ random point $z_i$, independently and uniformly distributed on the unit disc, we set:
\begin{eqnarray*}
T^{k,\beta}_m=\frac{1}{m}\sum_{i=1}^m \left(\sum_{j=1}^n \Delta f(z_i)(\log\lambda^k_j(z_i)-\log\gamma_j(z_i))\right)^\beta
\end{eqnarray*}

We have $||\Delta f||_{\infty}\leq C$ , and since logarithmic singularities are integrable:
\begin{eqnarray*}
\sup_{\beta=1,2}\int \left|\sum_{j=1}^n (\log\lambda^k_j(z)-\log\gamma_j(z))\right|^\beta\leq Cn^2
\end{eqnarray*}

Consequently, by the sampling lemma, for any $A >0$ we have with probability at least $1-n^{-A}$:

\begin{eqnarray*}
\sup_{\beta=1,2}\left|T^{k,\beta}_m-\int \Delta f(z) \left(\sum_{j=1}^n (\log\lambda^k_j(z)-\log\gamma_j(z))\right)^\beta\right|
\leq C\frac{n^2}{\sqrt{mn^{-A}}}
\end{eqnarray*}

Setting $m = n^{A+5}$, both difference goes to zero in the limit.

Now, fix $D$ and $\epsilon >0$. The rigidity estimate (Lemma \ref{p3}) fails with probability at most $n^{-A-5-D}$, so we may apply it to each summand in $S^k_m$ to obtain that, with probability at least $1-n^{-D}$, the following
estimate holds:

\begin{eqnarray*}
\sup_{\beta=1,2}\left|\int \Delta f(z)\left(\sum\log\lambda_j(z)-\sum\log \gamma_j(z)\right)^\beta dz\right|
= O(n^{\epsilon})
\end{eqnarray*}

Consequently, if $z$ is chosen uniformly randomly from the unit disc, the above estimates imply that for any choice of $\epsilon >0$, the following variance estimate holds with overwhelming probability:
\begin{eqnarray*}
\mathbf{Var}\left[\Delta f(z)\left(\log|\lambda^k_j(z)|-\log|\gamma_j(z)|\right)\right]\leq C_{\epsilon}n^{\epsilon}
\end{eqnarray*}

Notice that the variance appearing here is with respect to the selection of $z$, not the entries of the matrices; notice also that the classical position of the eigenvalues does not, of course, depend on $k$.

We will apply this variance estimate in a moment. For some natural number $m$ to be determined, let $z_1,....,z_M$ again be uniformly random, independently selected points in the unit disc, and again consider the empirical average $S^k_M$:
\begin{eqnarray*}
S^k_M=\frac{1}{M}\sum_{i=1}^M \left(\sum_{j=1}^n \Delta f(z_i)(\log\lambda^k_j(z_i)-\log\gamma_j(z_i))\right)
\end{eqnarray*}

By the triangle inequality:
\begin{eqnarray*}
\left|G\left(
\int \Delta f(z)\sum_{j=1}^n\left(\log\lambda^1_j(z)-\log\gamma_j(z)\right)dz
\right)
-
G\left(
\int \Delta f(z)\sum_{j=1}^n\left(\log\lambda^2_j(z)-\log\gamma_j(z)\right)dz
\right)\right|\\
\leq |G(S^1_M)-G(S^2_M)|+||G^{'}||_{\infty}\sum_{k=1}^2
\left|
\int \Delta f(z)\sum_{j=1}^n\left(\log\lambda^k_j(z)-\log\gamma_j(z)\right)dz
-S_M^k\right|
\end{eqnarray*}

Since we may assume that the variance estimate above occurs surely, another application of the sampling lemma implies:
\begin{eqnarray*}
\mathbf{P}\left[\left|\int \Delta f(z)\sum_{j=1}^n\left(\log\lambda^k_j(z)-\log\gamma_j(z)\right)dz-S_M^k \right|\geq \delta\right]
\leq \frac{1}{M\delta^2}\mathbf{Var}[S^k_M]\leq C_\epsilon\frac{n^{\epsilon}}{M\delta^2}
\end{eqnarray*}

Choosing $\epsilon = 1/8$ and $\delta= n^{-1/32}$ allows us to choose $M=O(n^{1/4})$ and have with probability $1-O(n^{-1/16})$ that:
\begin{eqnarray*}
\left|\left(\int \Delta f(z)\sum_{j=1}^n\log\lambda^k_j(z)dz-
\int \Delta f(z)\sum_{j=1}^n\log\gamma_j(z)dz\right)
-S^k_M \right|\leq O(n^{-1/32})
\end{eqnarray*}

This clearly goes to zero as $n$ becomes large, and $||G^{'}||_{\infty}$ is finite, so to prove our desired estimate it suffices to show that for some positive constants $C>0$ and $A>0$:

\begin{eqnarray*}
\Bigg|\mathbf{E}G\left(\frac{1}{M}\sum_{i=1}^M\sum_{j=1}^n \Delta f(z_i)(\log\lambda^1_j(z_i)-\log\gamma_j(z_i))\right)\\
-\mathbf{E}G\left(\frac{1}{M}\sum_{i=1}^M \sum_{j=1}^n \Delta f(z_i)(\log\lambda^2_j(z_i)-\log\gamma_j(z_i))\right)\Bigg|\\
\leq Cn^{-A}
\end{eqnarray*}

We can reduce the problem again. The fundamental theorem of calculus allows us to write the log-determinant in terms of the Stiejles transform:
\begin{eqnarray*}
\log|\det(W_n(z))|=\log|\det(W_n(z)-in^{100})|-n\mbox{Im}\int_0^{n^{100}}s(z,i\eta)d\eta\\
=100n\log(n)+O(n^{-10})-n\mbox{Im}\int_0^{n^{100}}s(z,i\eta)d\eta
\end{eqnarray*}

Since $G$ has a bounded first derivative, we may ignore the contribution of terms which go to zero asymptotically, and by centering (or translating $G$) we may ignore deterministic terms. Therefore it is sufficent to show:

\begin{eqnarray*}
\Bigg|\mathbf{E}G\left(\frac{n}{M}\sum_{i=1}^M\Delta f(z_i)\mbox{Im}\int_0^{n^{100}}s_1(z_i,\sqrt{-1}\eta)d\eta\right)\\
-\mathbf{E}G\left(\frac{n}{M}\sum_{i=1}^M\Delta f(z_i)\mbox{Im}\int_0^{n^{100}}s_2(z_i,\sqrt{-1}\eta)d\eta\right)\Bigg|\\
\leq Cn^{-A}
\end{eqnarray*}

To prove this inequality, we would like to bring in the machinery of resolvent swapping, but  first we need to make sure that the resolvent swapping lemma applies. If $|t|=O(n^{1/2-\epsilon})$ (which follows from our assumptions about
the atom distributions of $M_n^1$ and $M_n^2$)
, it follows from the proof of Theorem 24 in \cite{TV2}
that the condition $|t|\times||R_0||_{(\infty,1)}=o(\sqrt{n})$ is satisfied with probability $1-n^{-1/4-c}$ for sufficiently small $c > 0$. Since $M=O(n^{1/4})$, we may safely assume that this condition is satisfied surely for every summand in $S_M^k$, therefore the lemma applies.

Fix two natural number $a$ and $b$, both between 1 and $n$.
Let $M_n^0$ be the matrix which agrees with $M_n^1$ everywhere except the $(a,b)$ and $(b,a)$ positions, where it has zeros. Let $V$ be the elementary matrix with ones in the $(a,b)$-th positions (or position, if $a=b$).  Let $\xi^1_{a,b}$ denote the $(a,b)$-th entry of $M_n^1$
and let $\xi^2_{a,b}$ denote the $(a,b)$-th entry of $M_n^2$, and let $\tilde{M}^1_n=M_n^0+\xi^2_{a,b}V$. Notice that trivially $M_n^1=M_n^0+\xi_{a,b}^1V$.

By Lemma \ref{p2}:
\begin{eqnarray*}
\sum_{i=1}^M\Delta f(z_i)\mbox{Im}\int_0^{n^{100}}s_1(z_k,\sqrt{-1}\eta)d\eta
=
\sum_{i=1}^M
\Delta f(z_i)\mbox{Im}\int_0^{n^{100}}
\Bigg(s_0(z_i)+\sum_{j=1}^4 n^{-j/2}c_j(\eta,z_i)[\xi^1_{a,b}]^j\\+O(n^{-5/2}|\xi^1_{a,b}|^{k+1}n^{c_0})
\min(n^{c_0},\frac{1}{n\eta})\Bigg)
\end{eqnarray*}

And:
\begin{eqnarray*}
\sum_{i=1}^M\Delta f(z_i)\mbox{Im}\int_0^{n^{100}}\tilde{s}_1(z_i,\sqrt{-1}\eta)d\eta
=
\sum_{i=1}^M\Delta f(z_i)\mbox{Im}\int_0^{n^{100}}
\Bigg(s_0(z_i)+\sum_{j=1}^4 n^{-j/2}c_j(\eta,z_i)[\xi^2_{a,b}]^j\\+O(n^{-5/2}|\xi^2_{a,b}|^{k+1}n^{c_0})
\min(n^{c_0},\frac{1}{n\eta})\Bigg)
\end{eqnarray*}

The coefficients $c_j$ satisfy:
\begin{eqnarray*}
n\int_0^{n^{100}}c_j(i,\eta)d\eta=O(n^{c_0})
\end{eqnarray*}

Notice that $\mathbf{E}[\xi^1_{a,b}]^\beta=\mathbf{E}[\xi^2_{a,b}]^\beta$ for $\beta=1,2,3,4$. Notice also that a trivial integration gives:

\begin{eqnarray*}
\left|\int_0^{n^{100}}\min(n^{c_0},\frac{1}{n\eta})d\eta\right| \leq C
\end{eqnarray*}

Consequently, for a small constant $c$ and $\beta=1,2,3,4$:
\begin{eqnarray*}
\mathbf{E}\left(
\sum_{i=1}^M\Delta f(z_i)\mbox{Im}\int_0^{n^{100}}s_1(z_k,\sqrt{-1}\eta)d\eta
\right)^\beta
-
\mathbf{E}\left(
\sum_{i=1}^M\Delta f(z_i)\mbox{Im}\int_0^{n^{100}}\tilde{s}_1(z_i,\sqrt{-1}\eta)d\eta
\right)^\beta
\\
=O(n^{-5/2+c})
\end{eqnarray*}

By Taylor expansion:
\begin{eqnarray*}
\mathbf{E}\Big[G\left(
\sum_{i=1}^M\Delta f(z_i)\mbox{Im}\int_0^{n^{100}}s_1(z_k,\sqrt{-1}\eta)d\eta
\right)
-
G\left(
\sum_{i=1}^M\Delta f(z_i)\mbox{Im}\int_0^{n^{100}}\tilde{s}_1(z_i,\sqrt{-1}\eta)d\eta
\right)\Big]
\\
=O(n^{-5/2+c})
\end{eqnarray*}

Summing over every choice of $a$ and $b$, and applying the triangle inequality:
\begin{eqnarray*}
\mathbf{E}\Big[G\left(
\sum_{i=1}^M\Delta f(z_i)\mbox{Im}\int_0^{n^{100}}s_1(z_k,\sqrt{-1}\eta)d\eta
\right)
-
G\left(
\sum_{i=1}^M\Delta f(z_i)\mbox{Im}\int_0^{n^{100}}s_2(z_i,\sqrt{-1}\eta)d\eta
\right)\Big]
\\
=O(n^2n^{-5/2+c})=O(n^{-1/2+c})
\end{eqnarray*}

If we choose $c$ sufficiently small this quantity vanishes, finishing the proof $\P$.

\section{Central Limit Theorems for GinOE}

\subsection{Preliminary Calculations and Background}

Here, we gather basic information on independent entry random matrices whose entries have real Gaussian distributions. This is the real Ginbre ensemble, or GinOE.
Our main reference for this entire section will be to the 2009 paper of Borodin and Sinclair \cite{BS}, where explicit formulas for the $k$-point correlation functions of the GinOE ensembles are calculated. For clarity and consistency, we will follow the notation
and conventions of that paper throughout our discussion of GinOE.

Perhaps the most striking feature of the spectrum of GinOE matrices is that, unlike ensembles with complex entries such as GinUE, the spectrum splits into a real and complex portion. That is, away from the real line,
we expect no eigenvalues to lie in any fixed set with
two-dimensional Lebesgue measure zero (therefore the distribution of eigenvalues is said to be dominated by the Lebesgue measure away from the real line). However, we \textit{do} expect a certain percentage of
eigenvalues to located on the real life itself, so the eigenvalue distribution is \textit{not} dominated by two dimensional Lebesgue meansure on the real line. 

Notice that the limiting distribution of eigenvalues is still uniform on the unit disc, and therefore is dominated by the two dimensional Lebesgue measure everywhere,
so in the limit the contribution of real eigenvalues is not felt. Notice also that we now have two sources of randomness for these eigenvalues: both the number of real eigenvalues, and then the position of those eigenvalues.

We will consider test functions with support either entirely on or entirely away from the real line.  Since the matrix entries are real, the complex eigenvalues necessarily come in
conjugate pairs. Consequently, we will restrict ourselves to functions supported in the upper half plane. The $k$-point correlation functions of GinOE are known for matrices with even dimension, so we restrict ourselves
to matrices with dimension $2n$.

To write down the $k$-point correlation functions, we will need to employ the machinery of Pfaffians. One can show that the determinant of a skew-symmetric matrix can always be written as the square of a polynomial,
this polynomial is referred to as the Pfaffian. Pfaffians can also be defined in a more practical manner:

\begin{eqnarray*}
\mbox{Pf}(M)=\frac{1}{2^n n!}\sum_{\sigma\in S_{2n}}\mbox{sgn}(\sigma)\prod_{i=1}^{n} M_{\sigma(2i-1),\sigma(2i)}
\end{eqnarray*}

Here, $S_{2n}$ is the symmetric group on $2n$ elements and $\mbox{sgn}(\sigma)$ denotes the signature of the permuation.

Since the Pfaffian squared equals the determinant, the determinant specifies the Pfaffian up to a sign. Since all the Pfaffians we will see occur as probability densities, however, their signs will always be positive.

While Pfaffians are complicated in general, the formulas for the Pfaffians of small matrices are tractable:

\[ \mbox{Pf} \left( \begin{array}{cc}
0 & a \\
-a & 0 \end{array}\right)=a\]

And:

\[ \mbox{Pf} \left( \begin{array}{cccc}
0 & a & b & c\\
-a & 0 & d & e\\
-b & -d & 0 & f \\
-c & -e & -f & 0 \end{array}\right)=af-be+dc\]

By \cite{BS},  the $k$-point correlation functions of the real $2n$ dimensional Ginibre ensemble are given by:
\begin{eqnarray*}
p_k(x_1,...,x_k)=
\mbox{Pf}(K(x_i,x_j))_{1\leq i,j\leq k}
\end{eqnarray*}

We have defined the submatrices $K(x_i,x_j)$ as follows:
\begin{eqnarray*}
K(x_i,x_j)=\left( \begin{array}{cc}
D_{2n}(x_i,x_j) & S_{2n}(x_i,x_j) \\
-S_{2n}(x_j,x_i) & I_{2n}(x_i,x_j) \end{array} \right)
\end{eqnarray*}

The eigenvalues of GinOE are therefore said to constitute a\textit{ Pfaffian point process}. The entries of $K(x_i,x_j)$ depend on the locations of the arguments in the complex plane. We will look at two cases: complex/complex (both arguments have positive imaginary components) and real/real (both arguments have no imaginary component).

In the complex/complex case (which corresponds to $f$ being supported in the upper half plane, away from the real line), we have:

\begin{eqnarray*}
S_{2n}(z,z)=
\frac{i e^{-(1/2)(z-\bar{z})^2}}{\sqrt{2\pi}}(\bar{z}-z)
G(z,z)s_{2n}(z\bar{z})\\
S_{2n}(z,w)=\frac{i e^{-(1/2)(z-\bar{w})^2}}{\sqrt{2\pi}}(\bar{w}-z)
G(z,w)s_{2n}(z\bar{w})\\
D_{2n}(z,w)=
\frac{e^{-(1/2)(z-w)^2}}{\sqrt{2\pi}}(w-z)G(z,w)s_{2n}(zw)\\
I_{2n}(z,w)=\frac{ e^{-(1/2)(\bar{z}-\bar{w})^2}}{\sqrt{2\pi}}(\bar{z}-\bar{w})
G(z,w)s_{2n}(\bar{z}\bar{w})
\end{eqnarray*}

Here we have used the following definitions:

\begin{eqnarray*}
G(z,w)=\sqrt{\mbox{erfc}(\sqrt{2}\Im (z))\mbox{erfc}(\sqrt{2}\Im(w))}\\
s_{2n}(z)=e^{-z}\sum_{j=1}^{2n-2}\frac{z^{j}}{j!}
\end{eqnarray*}

This is all we need to compute an explicit representation of our linear statistic. Indeed:
\begin{eqnarray*}
\mbox{Var}[S_{2n}(f)]=\mathbf{E}\sum_{j=1}^{2n} f\left(\frac{\lambda_j}{\sqrt{2n}}\right)^2
+\mathbf{E}\sum_{j\neq k}
f\left(\frac{\lambda_j}{\sqrt{2n}}\right)
f\left(\frac{\lambda_k}{\sqrt{2n}}\right)
-\left(\mathbf{E}\sum_{j=1}^{2n} f\left(\frac{\lambda_j}{\sqrt{2n}}\right)\right)^2
\end{eqnarray*}

Or, in terms of correlation functions:
\begin{eqnarray*}
\mbox{Var}(S_{2n}(f))=
\int f\left(\frac{x}{\sqrt{2n}}\right)^2 p^{(2n)}_1(x)dx
+\int f\left(\frac{x}{\sqrt{2n}}\right)f\left(\frac{y}{\sqrt{2n}}\right)p^{(2n)}_2(x,y)\\
-\left(\int f\left(\frac{x}{\sqrt{2n}}\right) p^{(2n)}(x)dx\right)^2
\end{eqnarray*}

After substituting in our explicit formulas, and after some small amount of cancellation and substitution, we have:

\begin{eqnarray*}
\mbox{Var}(S_{2n}(f))=
2n\int f(x)^2S_{2n}(\sqrt{2n}x,\sqrt{2n}x)dx-
4n^2\int f(x)f(y)\\
\times\left(
D_{2n}(\sqrt{2n}x,\sqrt{2n}y)I_{2n}(\sqrt{2n}x,\sqrt{2n}y)
+S_{2n}(\sqrt{2n}x,\sqrt{2n}y)S_{2n}(\sqrt{2n}y,\sqrt{2n}x)
\right)dxdy
\end{eqnarray*}

In the real/real case, where $f$ is supported on the real line, we have:

\begin{eqnarray*}
S_{2n}(x,y)=\frac{1}{\sqrt{2\pi}}e^{-\frac{1}{2}(x-y)^2}
e^{-xy}\sum_{m=0}^{2n-2}\frac{(xy)^m}{m!}
\\
+e^{-x^2/2}\frac{2^{n-3/2}}{\sqrt{\pi}(2n-2)!}\mbox{sgn}(x)x^{2n-1}\gamma(n-\frac{1}{2},y^2/2)
\end{eqnarray*}

Here, $\gamma(x,y)$ denotes the lower incomplete gamma function:
\begin{eqnarray*}
\gamma(x,y)=\int_0^{y} t^{x-1}e^{-x}dx
\end{eqnarray*}

We also define:
\begin{eqnarray*}
D_{2n}(x,y)=(y-x)\frac{1}{\sqrt{2\pi}}e^{-\frac{1}{2}(x-y)^2}e^{-xy}\sum_{m=0}^{2n-2}\frac{(xy)^m}{m!}
\end{eqnarray*}

And finally:
\begin{eqnarray*}
I_{2n}(x,y)=e^{-x^2/2}\frac{1}{2\sqrt{\pi}}\int_0^{y^2/2}\frac{e^{-t}}{\sqrt{t}}
\sum_{m=0}^{n-1}\frac{(\sqrt{2t}x)^{2m}}{(2m)!}dt\\-
e^{-y^2/2}\frac{1}{2\sqrt{\pi}}\int_0^{x^2/2}\frac{e^{-t}}{\sqrt{t}}
\sum_{m=0}^{n-1}\frac{(y\sqrt{2t})^{2m}}{(2m)!} dt
+\frac{1}{2}\mbox{sgn}(x-y)
\end{eqnarray*}

Plugging in the one and two point correlation functions, we obtain the following expression for the variance of the linear statistic:

\begin{eqnarray*}
\mbox{Var}(S_{2n}(f))=
\sqrt{2n}\int f(x)^2S_{2n}(\sqrt{2n}x,\sqrt{2n}x)dx-
2n\int f(x)f(y)\\
\times\left(
D_{2n}(\sqrt{2n}x,\sqrt{2n}y)I_{2n}(\sqrt{2n}x,\sqrt{2n}y)
+S_{2n}(\sqrt{2n}x,\sqrt{2n}y)S_{2n}(\sqrt{2n}y,\sqrt{2n}x)
\right)dxdy
\end{eqnarray*}

We will deal with these variance formulas shortly, but first, we need to fix some notation regarding Fourier transforms. We will find it convenient to define Fourier transform to be:
\begin{eqnarray*}
F[f](\omega)=\int f(x)e^{-i\omega x}dx
\end{eqnarray*}

In other words, we will write things in terms of angular frequency. The inverse transform is then:
\begin{eqnarray*}
f(x)=\frac{1}{2\pi}\int F[f](\omega)e^{ix\omega}d\omega
\end{eqnarray*}

Under this definition, we have the identities:
\begin{eqnarray*}
F[e^{-ax^2}](\omega)=\sqrt{\frac{\pi}{a}}e^{-k^2/4a}\\
F[f(x)g(x)]=\frac{1}{2\pi}\int f(\omega-\tau)g(\tau)d\tau\\
F[f*g]=F[f](\xi)\times F[g](\xi)\\
\int f(x)g(x)dz= \frac{1}{2\pi}\int \overline{\hat{f}(x)}\hat{g}(x)dx\\
i\xi \hat{f}(\xi)=F\left[\frac{d}{dx}f(x)\right]
\end{eqnarray*}

Here, $f$ and $g$ are Schwartz functions.

\subsection{A Detour For Quaternions}

Here we develop elements from the theory of quaternions which we will need later. We will follow the approach and notation of \cite{Meh}.

Define the algebra of \textit{complexified quaternions} to be the four dimensional vector space over the field of complex
numbers with basis $1, e_1, e_2, e_3$, equipped with the operation of quaternion multiplication:

\begin{eqnarray*}
e_1^2= e_2^2=e_3^2=e_1e_2e_3=-1
\end{eqnarray*}

This algebra is non-commutative, and is easily seen to be isomorphic to the algebra of two by two matrices by making the following indentifications:

\[ e_1= \left( \begin{array}{cc}
0 & i \\
i & 0 \end{array}\right)\]

\[ e_2 \left( \begin{array}{cc}
0 & -1 \\
1 & 0 \end{array}\right)\]

\[ e_3= \left( \begin{array}{cc}
i & 0 \\
0 & -i \end{array}\right)\]

We can use this identification to define conjugates of quaternions. The relationship between a generic quaternion $q$ and its conjugate quaternion $q^{*}$  is given as follows:

\[ q= \left( \begin{array}{cc}
a & b \\
c & d \end{array}\right)\]

\[ q^{*}= \left( \begin{array}{cc}
d & -b \\
-c & a \end{array}\right)\]

We will consider matrices with complexified quaternion elements. By the identification above, one sees that an $n$ by $n$ quaternion matrix $Q$ corresponds to a $2n$ by $2n$ complex matrix. We let $\phi(Q)$ denote
this second matrix.

If $Q$ is a quaternion matrix, then a right eigenvalue of $Q$ is
defined to be a scalar $\lambda$ such that $Qv=\lambda v$ for some non-zero quaternion vector $v$. If the entries of $Q$ satisfy
$q_{ij}=q_{ji}^{*}$ for all $i$ and $j$, then $Q$ is said to be self-dual and has exactly $n$ eigenvalues, which are all real.

There are several ways to define the determinant of a quaternion matrix $Q$. For instance, one way is to simply swap quaternions with their matrix representations and take the usual determinant of the
resulting $2n$ by $2n$ matrix $\phi(Q)$-- this is known as the Study determinant. Another option is the Moore-Dyson determinant, which is defined to be the product of eigenvalues. This is the determinant which
we will find useful, and we will denote it $\det (Q)$. 

Notice that the Study determinant and the Moore-Dyson determinant do not, in general, agree with one another. In particular, while the Study determinant is a complex number by definition, the Moore Dyson determinant may be a quaternion.

To use $\det (Q)$ effectively, we will need a way to compute this quantity without having to find any eigenvalues. Consider some permutation $\sigma$ which is the product $m$ of disjoint cycles:

\begin{eqnarray*}
\sigma=(n_1k^1_2...k^1_{s_1})...(n_rk^r_2...k^r_{s_m})
\end{eqnarray*}

We will henceforth adopt the convention that $n_i$ is the largest entry in each cycle, and $n_i\geq n_{i+1}$.

Let $S_r$ be the symmetric group on $r$ elements. As an analogue to the classical Cayley combinatorial formula for determinants,
we have the following equivalent definition of the Moore-Dyson determinant of the quaternion matrix $Q$:

\begin{eqnarray*}
\det (Q)=\sum_{\sigma\in S_r}(-1)^{\mbox{sgn}(\sigma)}(q_{n_1i_2}...q_{i_sn_1})...
(q_{n_rk_2}...q_{k_hn_m})
\end{eqnarray*}

Let's assume now that $Q$ is a self-dual quaternion matrix, this gives us at least two more useful properties. First, the Moore Dyson determinant is necessarily complex valued for self-dual matrices. For the second, let $Z$ is the matrix with blocks $\left( \begin{array}{cc}
0 & 1 \\
-1 & 0 \end{array}\right)$ on the main diagonal and zeros elsewhere (thus $Z$ is the tensor product of the identity and this two dimensional matrix). Then $Z\phi(Q)$ is skew-symmetric, and we have
the following formula due to Mehta:

\begin{eqnarray*}
\det(Q)=\mbox{Pf}(Z \phi(Q))
\end{eqnarray*}

This relationship allows us to rewrite Pfaffian formulas in terms of determinantal formulas, which in turn will allow us to employ the powerful theory available for
determinantal point processes. To apply these, we will need to introduce cumulants.

The $N$-th cumulant of a random variable $X$, written $C_N$, is defined by the formula:

\begin{eqnarray*}
\log \mathbf{E}e^{tX}=\sum_{n=1}^{\infty}C_{N}\frac{t^n}{n!}
\end{eqnarray*}

Cumulants uniquely determine a random variable (by way of the Fourier inversion formula), and it can be shown that all cumulants of a normal random variable after the second are identically zero. 

We will be able to compute cumulants via the following lemma:

\begin{lemma}
\label{nh:sos}
Let  $Q(\cdot,\cdot)$ be a function which takes values in the algebra of complexified quaternions.
Let $(\lambda_1,...,\lambda_n)$ be a random set of $n$ points whose correlation functions satisfy:

\begin{eqnarray*}
p_r(x_1,...,x_r)=\det Q_r=\det \left(Q(x_a,x_b)\right)_{1\leq a,b\leq r}
\end{eqnarray*}

Then, if $f$ is a test function and $S[f]=\sum_{i=1}^n f(\lambda_i)$ is the linear statistic, the $k$-th cumulant of $S[f]$ is given by:
\begin{eqnarray*}
C_{k}(S[f])=\int\Bigg(\sum_{V_i\in X_k}
f(x_{1})^{|V_{i,1}|}...f(x_{|V_i|})^{|V_{i,|V_i|}|}
(-1)^{|V_i|-1}\\
\times
\sum_{\sigma\in C[|V_i|]}
\left(Q_n(x_{1},x_{\sigma(1)})Q_n(x_{\sigma(1)},x_{\sigma^2(1)})...
Q_n(x_{\sigma^{-1}(1),x_1})
\right)
\Bigg)dx_1...dx_{|V_i|}
\end{eqnarray*}

 Here, $C[r]$ denotes the collection of circular permutations on $r$ elements, and $X^k$ denotes the collection of paritions of the set $[1,...,k]$.
\end{lemma}

The proof is identical to the proof of the analagous statement for complex-valued functions, which is a classical result of Soshnikov \cite{So} (commutivity does not come up in the proof).

\subsection{Variance in the Complex/Complex Case}

In this section, we will compute the variance of the linear statistic in the complex/complex case. In particular, we will show:
\begin{proposition}
We have the following asymptotic formula:
\begin{eqnarray*}
\lim_{n\to\infty}\text{\emph{Var}}(S_{2n}(f))=\frac{1}{4\pi}\int (f_x^2(z)+f_y^2(z)) dz
\end{eqnarray*}
Here $f$ is a Schwartz function with compact support in the upper half plane, supported away from the real line.
\end{proposition}

We will prove this by a series of lemmas. Repeatedly, we will make use of the
following lemma, found in \cite{BS}:

\begin{lemma}
Let $u$ be a complex number in a compact subset of the open unit disc. Then:
\begin{eqnarray*}
s_{2n}(2nu)=1-\frac{1}{\sqrt{n}}
\frac{e^{-2(1-u)}}{2\pi u }e^{2n(1-u)}u^{2n}\left(1+O\left(\frac{1}{n}\right)\right)
\end{eqnarray*}
In particular, if $u$ is real and $0<u<1$, then $s_{2n}(2nu)=1+O(1/n)$.
\end{lemma}

In practice, this will allow us to replace $s_{2n}(2nu)$ with 1 during our computations.

The first term in the variance expression is handled in the following easy lemma:

\begin{lemma}
We have the following identity (in the complex/complex case):
\begin{eqnarray*}
2n\int f(x)^2S_{2n}(\sqrt{2n}x,\sqrt{2n}x)dx=\frac{2n}{\pi}
\int f(z)^2dz-\frac{1}{4\pi}\int \frac{f(z)^2}{|\Im(z)|^2}dz+O(n^{-1})
\end{eqnarray*}
Here, we take $f$ to be compactly supported inside the upper half of the open unit disc, away from the real line,
and square integrable.
\end{lemma}

\textbf{Proof:}

We have the well-known asymptotic for large $x$:
\begin{eqnarray*}
\mbox{erfc}(x)=1-\mbox{erf}(x)=\frac{e^{-x^2}}{\sqrt{\pi}x}\left(1-\frac{1}{2x^2}+O(\frac{1}{x^4})\right)
\end{eqnarray*}

Consequently:
\begin{eqnarray*}
G(\sqrt{2n}z,\sqrt{2n}z)=\sqrt{\mbox{erf}(\sqrt{2}\sqrt{2n}\Im(z),\sqrt{2}\sqrt{2n}\Im(z))^2}
\\=\frac{e^{-4n\Im(z)^2}}{\sqrt{2n}\sqrt{2\pi}|\Im(z)|}\left(1-\frac{1}{8n|\Im(z)|^2}+O(n^{-2})\right)
\end{eqnarray*}

We have:
\begin{eqnarray*}
2n\int f(z)^2S_{2n}(\sqrt{2n}z,\sqrt{2n}z)dz=
2n\sqrt{2n}\int f(z)^2
\frac{i e^{-4n(i\Im(z))^2}}{\sqrt{2\pi}}(-2i\Im(z))
G(\sqrt{2n}z,\sqrt{2n}z)s_n(2nz\bar{z})dz\\
=4n\sqrt{2n}
\int f(z)^2
\frac{ e^{4n\Im(z)^2}}{\sqrt{2\pi}}(\Im(z))
G(\sqrt{2n}z,\sqrt{2n}z)s_{2n}(2n|z|)dz\\
=
\frac{2n}{\pi}
\int f(z)^2\left(1-\frac{1}{8n|\Im(z)|^2}+O(n^{-2})\right)
s_{2n}(2n|z|)dz
\end{eqnarray*}

Since $s_{2n}(2n|z|)=1+O(e^{-cn})$, the proof is completed $\P$.

With one term taken care of, we only have two more to deal with.

\begin{lemma}
We have the identity (in the complex/complex case):

\begin{eqnarray*}
4n^2\int f(x)f(y)
D_{2n}(\sqrt{2n}x,\sqrt{2n}y)I_{2n}(\sqrt{2n}x,\sqrt{2n}y)dxdy=O(e^{-nc})
\end{eqnarray*}

We take $f$ to be compactly supported in the upper half of the open unit disc, away from the real line, and bounded.
\end{lemma}

\textbf{Proof:}

We will show that both $D_{2n}(\sqrt{2n}z,\sqrt{2n}w)$ and $I_{2n}(\sqrt{2n}w,\sqrt{2n}w)$ vanish exponentially.
We have the definitions:

\begin{eqnarray*}
D_{2n}(\sqrt{2n}z,\sqrt{2n}w)=\sqrt{2n}
\frac{e^{-n(z-w)^2}}{\sqrt{2\pi}}(w-z)G(\sqrt{2n}z,\sqrt{2n}w)s_{2n}(2nzw)\\
I_{2n}(\sqrt{2n}w,\sqrt{2n}w)=\sqrt{2n}\frac{e^{-n(\bar{z}-\bar{w})^2}}{\sqrt{2\pi}}(\bar{z}-\bar{w})
G(\sqrt{2n}z,\sqrt{2n}w)s_{2n}(2n\bar{z}\bar{w})
\end{eqnarray*}

We will show that $D_{2n}(\sqrt{2n}z,\sqrt{2n}w)$ decays uniformly exponentially
on the support of $f$, the proof for $I_{2n}(\sqrt{2n}w,\sqrt{2n}w)$ being essentially identical (only the overall phase will change,
which does not impact the proof). We begin with the trivial expansion:
\begin{eqnarray*}
e^{-n(z-w)^2}=e^{-n((\Re(z)-\Re(w))^2-(\Im(z)-\Im(w))^2)}
\end{eqnarray*}

We may also expand:
\begin{eqnarray*}
G(\sqrt{2n}z,\sqrt{2n}w)= \left(
\frac{e^{-4n\Im(z)^2}e^{-4n\Im(w)^2}}{4n\pi\Im(z)^2\Im(w)^2}
\right)^{1/2}\left(1+O\left(\frac{1}{n}\right)\right)\\
=\frac{e^{-2n(\Im(z)^2+\Im(w)^2)}}{2\sqrt{n\pi}|\Im(z)\Im(w)|}\left(1+O\left(\frac{1}{n}\right)\right)
\end{eqnarray*}

And:
\begin{eqnarray*}
s_{2n}(2nzw)= \left(1-O(1)\frac{1}{\sqrt{n}}e^{2n(1-zw)}zw^{2n}\right)\left(1+O\left(\frac{1}{n}\right)\right)
\end{eqnarray*}

Consequently:
\begin{eqnarray*}
\left|D_{2n}(\sqrt{2n}z,\sqrt{2n}w)\right|\leq C
e^{-n((\Re(z)-\Re(w))^2-(\Im(z)-\Im(w))^2)}\\
e^{-2n(\Im(z)^2+\Im(w)^2)}
 \left(1-O(1)\frac{1}{\sqrt{n}}e^{2n(1-zw)}zw^{2n}\right)
\left(1+O\left(\frac{1}{n}\right)\right)
\end{eqnarray*}

We have:
\begin{eqnarray*}
e^{-n((\Re(z)-\Re(w))^2-(\Im(z)-\Im(w))^2)}
e^{-2n(\Im(z)^2+\Im(w)^2)}\\
=e^{-n((\Re(z)-\Re(w))^2+\Im(z)^2+\Im(w)^2+2\Im(z)\Im(w))}\\
\leq C e^{-2n\Im(z)\Im(w)}\leq Ce^{-2nc}
\end{eqnarray*}

We now consider the term:
\begin{eqnarray*}
e^{-n((\Re(z)-\Re(w))^2-(\Im(z)-\Im(w))^2)}e^{2n(1-zw)}zw^{2n}
e^{-2n(\Im(z)^2+\Im(w)^2)}
\end{eqnarray*}

Combining exponential terms, we get an oscillating term of magnitude one, multiplied by an amplitude of:
\begin{eqnarray*}
\mbox{exp}\Big[-n\Big(
(\Re(z)-\Re(w))^2-(\Im(z)-\Im(w))^2
-2(1-\Re(zw)+\log|zw|)\\
+2\Im(z)^2+2\Im(w)^2
\Big)
\Big]
\end{eqnarray*}

Or else:
\begin{eqnarray*}
\mbox{exp}\Big[-n\Big(
(\Re(z)-\Re(w))^2-(\Im(z)-\Im(w))^2
-2(1+\log|zw|)\\
+2(\Re(z)\Re(w)-\Im(z)\Im(w))+2\Im(z)^2+2\Im(w)^2
\Big)
\Big]\\
=
\mbox{exp}\Big[-n\Big(
\Re(z)^2+\Re(w)^2
-2(1+\log|zw|)+\Im(z)^2+\Im(w)^2
\Big)
\Big]
\end{eqnarray*}

We may rewrite this last display:
\begin{eqnarray*}
\mbox{exp}\Big[-n(\Big( (|z|^2-\log|z|^2-1)+(|w|^2-\log|w|^2-1)\Big) \Big]
\end{eqnarray*}

This last quantity in the exponent is necessarily positive (because $|z|$ and $|w|$ are less than one by assumption),
 and consequently vanishes exponentially, proving the result. $\P$

It remains to control the last term, which we turn to presently.

\begin{lemma}
We have the following expansion (in the complex/complex case):
\begin{eqnarray*}
4n^2\int f(x)f(y)
S_{2n}(\sqrt{2n}x,\sqrt{2n}y)S_{2n}(\sqrt{2n}y,\sqrt{2n}x)dxdy
\\
=\frac{2n}{\pi}\int f(z)^2 dz-\frac{1}{4\pi}\int (f_x^2(z)+f_y^2(z)) dz
-\frac{1}{4\pi}\int \frac{f(z)^2}{\Im(z)^2}dz+o(1)
\end{eqnarray*}
Here, $f$ is supported in the upper half of the open unit disc, away from the real line, and is a Schwartz function.
\end{lemma}

\textbf{Proof:}

By definition, we have:
\begin{eqnarray*}
4n^2\int\int f(z)f(w)S_{2n}(\sqrt{2n}x,\sqrt{2n}y)S_{2n}(\sqrt{2n}y,\sqrt{2n}x)dxdy\\
=8n^3\int\int f(z)f(w)\frac{e^{-n[(z-\bar{w})^2+(w-\bar{z})^2]}}{2\pi}|\bar{w}-z|^2
G(\sqrt{2n}z,\sqrt{2n}w)^2 s_n(2nz\bar{w})s_n(2nw\bar{z})
\end{eqnarray*}

If we expand $G$ using the usual asymptotic estimate for $\mbox{erfc}(x)$ and replace $s_n(2nz\bar{w})s_n(2nw\bar{z})$ with one, the magnitude
of the error term incurred is of the order:
\begin{eqnarray*}
\mbox{exp}\Big[-2n\Big(|z-w|^2-2+2\Re(z\bar{w})-2\log(|zw|)\Big)\Big]\\
=
\mbox{exp}\Big[-2n\Big(|z|^2+|w|^2-2-\log|z|^2-\log|w|^2\Big)\Big]
\end{eqnarray*}

Since $|z|$ and $|w|$ are less than one, this vanishes exponentially. We are left with the following expression at leading order:

\begin{eqnarray*}
8n^3\int\int f(z)f(w)
\frac{e^{-2n\Re[(z-\bar{w})^2]}}{2\pi}|\bar{w}-z|^2
\frac{e^{-4n\Im(z)^2}e^{-4n\Im(w)^2}}{4n\pi |\Im(z)|\times|\Im(w)|}\\
\left(1-\frac{1}{8n|\Im(z)|^2}+O(n^{-2})\right) \left(1-\frac{1}{8|\Im(w)|^2}+O(n^{-2})\right)
dzdw
\end{eqnarray*}

Remembering $\Im(z),\Im(w)>0$, we can rewrite:
\begin{eqnarray*}
\frac{|\bar{w}-z|^2}{\Im(z)\Im(w)}=4+\frac{|w-z|^2}{\Im(z)\Im(w)}\\
\Re[(z-\bar{w})^2]+2\Im(z)^2+2\Im(w)^2=|z-w|^2
\end{eqnarray*}

Our leading order expression is now:

\begin{eqnarray*}
\frac{n^2}{\pi^2}\int\int f(z)f(w)
e^{-2n|z-w|^2}\left(4+\frac{|w-z|^2}{\Im(z)\Im(w)} \right)\\
\left(1-\frac{1}{8n|\Im(z)|^2}\right) \left(1-\frac{1}{8n|\Im(w)|^2})\right)
dzdw
\end{eqnarray*}

To control this integral, we will prove the following pair of identities:
\begin{eqnarray*}
\int\int f(z)g(w) e^{-2n|z-w|^2}dzdw=
\frac{\pi}{2n}\int\int f(x,y)g(x,y)dxdy\\
-\frac{\pi}{16n^2}\int\int \left(\frac{df}{dx}\frac{dg}{dx}+ \frac{df}{dy}\frac{dg}{dy}\right)dxdy+O(n^{-3})
\end{eqnarray*}

And:
\begin{eqnarray*}
\frac{4n^2}{\pi}\int f(z)\int f(w)|w-z|^2 e^{-2n|w-z|^2}dwdz=\int f(z)^2 dz+o(1)
\end{eqnarray*}

Substituting these into our leading order expression, we obtain, as desired:
\begin{eqnarray*}
\frac{2n}{\pi}\int f(z)^2 dz-\frac{1}{4\pi}\int (f_x^2(z)+f_y^2(z)) dz
-\frac{1}{4\pi}\int \frac{f(z)^2}{\Im(z)^2}dz
\end{eqnarray*}

It remains to prove the two identities claimed. The first is an elementary Fourier transform argument,
while the second is an approximation to the identity. Indeed:

\begin{eqnarray*}
\int\int f(z)g(w) e^{-2n|z-w|^2}dzdw=\frac{1}{4\pi^2}\int \overline{F[f(u)](\xi)}
F[g(u)*e^{-2nu^2}](\xi)d\xi\\
=\frac{1}{4\pi^2}\frac{\pi}{2n}\int \overline{\hat{f}(\xi)}\hat{g}(\xi)e^{-\xi^2/8n}d\xi
=\frac{1}{4\pi^2}\frac{\pi}{2n}\int \overline{\hat{f}(\xi)}\hat{g}(\xi)\left(1-\frac{\xi^2}{8n}+O(n^{-2})\right)d\xi
\end{eqnarray*}

Or else:
\begin{eqnarray*}
\frac{\pi}{2n}\int \Big[f(z)g(z)-\frac{d}{dz}f(z)\times \frac{d}{d\bar{z}}g(z)\frac{1}{8n}+O(n^{-2})\Big]d^2z
\end{eqnarray*}

To prove the second claim, start with:

\begin{eqnarray*}
n^2\int_{\mathbf{C}} |z|^2 e^{-2n|z|^2}dz=2\pi n^2 \int_0^{\infty} r^3 e^{-2nr^2}dr\\
=2 \pi \int_0^\infty \rho^3 e^{-2\rho ^2} d\rho=\frac{2\pi}{8}=\frac{\pi}{4}
\end{eqnarray*}

Consequently, $\frac{4n^2}{\pi}|z|^2e^{-2n|z|^2}$ is an approximation of the identity, and so:
\begin{eqnarray*}
\frac{4n^2}{\pi}\int f(z)\int f(w)|w-z|^2 e^{-2n|w-z|^2}dwdz=\int f(z)^2 dz+o(1)
\end{eqnarray*}

The $o(1)$ term vanishes as $n$ becomes large, so this is what we wanted to show $\P$.

Substituting these lemmas into our expression for the variance, we have:

\begin{eqnarray*}
\mbox{Var}(S_{2n}(f))=
2n\int f(x)^2S_{2n}(\sqrt{2n}x,\sqrt{2n}x)dx-
4n^2\int f(x)f(y)\\
\times\left(
D_{2n}(\sqrt{2n}x,\sqrt{2n}y)I_{2n}(\sqrt{2n}x,\sqrt{2n}y)
+S_{2n}(\sqrt{2n}x,\sqrt{2n}y)S_{2n}(\sqrt{2n}y,\sqrt{2n}x)
\right)dxdy\\
=\frac{2n}{\pi}
\int f(z)^2dz-\frac{1}{4\pi}\int \frac{f(z)^2}{|\Im(z)|^2}dz
-\Big(\frac{2n}{\pi}\int f(z)^2 dz-\frac{1}{4\pi}\int (f_x^2(z)+f_y^2(z)) dz\\
-\frac{1}{4\pi}\int \frac{f(z)^2}{\Im(z)^2}dz\Big)+O(n^{-1})
\end{eqnarray*}

And of course this last display is just:
\begin{eqnarray*}
\frac{1}{4\pi}\int (f_x^2(z)+f_y^2(z)) dz+O(n^{-1})
\end{eqnarray*}

This completes the proof of the proposition.

\subsection{Higher Cumulants in the Complex/Complex Case: Proof of Theorem \ref{nh:bulk}}

Now that we have computed the limiting variance, the next step is to show that the limiting distribution of the linear statistic is a Gaussian random variable.
To do this, it is sufficient to show that all cumulants after the second vanish in the limit.

We will begin with some simplifications. Recall that the $k$-point correlation functions of GinOE, away from the real line, take the form of a Pfaffian of a matrix
constructed out of the following submatrices:

\begin{eqnarray*}
K_{2n}(x_i,x_j)=\left( \begin{array}{cc}
D_{2n}(x_i,x_j) & S_{2n}(x_i,x_j) \\
-S_{2n}(x_j,x_i) & I_{2n}(x_i,x_j) \end{array} \right)
\end{eqnarray*}

In the previous section, we have shown that $D_{2n}(x_i,x_j)$ and $I_{2n}(x_i,x_j)$ decay exponentially in $n$. Since the Pfaffian is a polynomial, and since
no terms grow exponentially, we instead consider the Pfaffian point process whose correlation functions are defined as follows:

\begin{eqnarray*}
\tilde{p}_{2n,k}(x_1,...,x_k)=
\mbox{Pf}(\tilde{K}_{2n}(x_i,x_j))_{1\leq i,j\leq k}
\end{eqnarray*}

The submatrices above are defined:

\begin{eqnarray*}
\tilde{K}_{2n}(x_i,x_j)=\left( \begin{array}{cc}
0 & S_{2n}(x_i,x_j) \\
-S_{2n}(x_j,x_i) & 0 \end{array} \right)
\end{eqnarray*}

We will denote the linear statisic associated with this point process by $\Lambda_{2n}[f]$.  By Lemma \ref{nh:sos}, if we can prove that all higher cumulants vanish asymptotically for this
random variable, we also obtain the result for $S_{2n}(f)$. To this end, consider the following matrix product:

\begin{eqnarray*}
Z\tilde{K}_{2n}(x_i,x_j)=\left( \begin{array}{cc}
S_{2n}(x_j,x_i) & 0 \\
0 & -S_{2n}(x_i,x_j) \end{array} \right)
\end{eqnarray*}

The quaternion matrix associated to $Z\tilde{K}_{2n}$, denoted $Q[Z\tilde{K}_{2n}]$, is self dual, since the functions $S_{2n}(x_i,x_j)$ satisfy:
\begin{eqnarray*}
S_{2n}(z_i,z_j)=-S_{2n}(z_j,z_i)
\end{eqnarray*}

Using $Z^2=I$, we may write things in terms of the Moore-Dyson determinant:
\begin{eqnarray*}
\mbox{Pf}(\tilde{K}_{2n})=\mbox{Pf}[Z^2 \tilde{K}_{2n}]=\det[Z\tilde{K_{2n}}]
=\det (S_n(x_i,x_j))_{1\leq i,j \leq n}
\end{eqnarray*}

Applying Lemma \ref{nh:sos}, we obtain an expression for the $k$-th cumulant of the random variable $\Lambda_{2n}[f]$:

\begin{eqnarray*}
\mathcal{C}_k(\Lambda_{2n}[f])=\int\Bigg(\sum_{V_i\in X_k}
f(z_{1})^{|V_{i,1}|}...f(z_{|V_i|})^{|V_{i,|V_i|}|}
(-1)^{|V_i|-1}\\
\times
\sum_{\sigma\in C[|V_i|]}
\left(S_n(z_{1},z_{\sigma(1)})S_n(z_{\sigma(1)},z_{\sigma^2(1)})...
S_n(z_{\sigma^{-1}(1)},z_1)
\right)
\Bigg)d^1(z_1)...d^2(z_{|V_i|})
\end{eqnarray*}

We can rearrange this as a sum over the size of partitions:

\begin{eqnarray*}
\sum_{m=1}^{k} (-1)^{m-1} \int\Bigg(\sum_{V_i\in X_N:|V_i|=m}
f(z_{1})^{|V_{i,1}|}...f(z_{m})^{|V_{i,m}|}\\
\times
\sum_{\sigma\in C[m]}
\left(S_n(z_{1},z_{\sigma(1)})S_n(x_{\sigma(1)},z_{\sigma^2(1)})...
S_n(z_{\sigma^{-1}(1)},z_1)
\right)
\Bigg)d^2(z_1)...d^2(z_m) 
\end{eqnarray*}

The well known \textit{multinomial coefficients} count the number of ways to put $N$ distinct objects into boxes of size $k_1$,...,$k_m$ with $k_1+....+ k_m=N$, which turns out to be:
\begin{eqnarray*}
\frac{N!}{k_1!...k_m!}
\end{eqnarray*}

Therefore (picking up a a factor of $1/m!$ due to the arbitrariness of partition labeling and a factor of $(m-1)!$ due to the sum over circular permutations) we can write:

\begin{eqnarray*}
\mathcal{C}_k(\Lambda_{2n}[f])=\sum_{m=1}^k\frac{(-1)^{m-1}}{m}\sum_{k_1+...+k_m=k}\frac{k!}{k_1!...k_m!}\\
\times
\int_{C^k}\left(\prod_{h=1}^m f(z_h)^{k_h}\right)S_{2n}(z_1,z_2)...S_{2n}(z_m,z_1)d^2(z_1)...d^2(z_m)
\end{eqnarray*}

This expression is sometimes known as the Costin-Lebowtiz formula for the cumulants of a determinantal point process. We are now in a position to prove:

\begin{proposition}\label{highc}
Let $f$ be a smooth function compactly supported in the upper half of the open unit disc, away from the real line. Then quantities
$\mathcal{C}_k(\Lambda_{2n}[f])$ vanish as $n$ becomes large for $k\geq 3$.
\end{proposition}

This implies a normal limiting distribution (as this is the only distribution for which all higher cumulants vanish).

The proof of the proposition is a combinatorial argument which will take the form of a series of lemmas. First, we have
a pair of simple combinatorial identities:

\begin{lemma}\label{com0}
For all integers $k>2$, the following identities hold:
\begin{eqnarray*}
\sum_{m=1}^k\frac{(-1)^{m-1}}{m}\sum_{k_1+...+k_m=k}\frac{k!}{k_1!...k_m!}=0\\
\sum_{m=1}^k\frac{(-1)^{m-1}}{m}\sum_{k_1+...+k_m=k}\frac{k!}{k_1!...k_m!}
\sum_{i \neq j}k_ik_j=0
\end{eqnarray*}
\end{lemma}

\textbf{Proof:}
We will use the method of generating functions. Set:

\begin{eqnarray*}
a_N=\sum_{m=1}^N\frac{(-1)^{m-1}}{m}\sum_{k_1+...+k_m=N}\frac{N!}{k_1!...k_m!}
\end{eqnarray*}

Then:
\begin{eqnarray*}
\Psi(x)=\sum_{N=1}^{\infty} a_N \frac{x^N}{N!}=\sum_{N=1}^{\infty}
\sum_{m=1}^N\frac{(-1)^{m-1}}{m}\sum_{k_1+...+k_m=N}\frac{x^N}{k_1!...k_m!}\\
=\sum_{m=1}^{\infty}\frac{(-1)^{m-1}}{m}\sum_{k_1>0, k_2 >0,...k_m>0}\frac{x^{k_1}...x^{k_m}}{k_1!...k_m!}
\end{eqnarray*}

Assuming $|e^{x}-1|\leq 1$ and taking the inner sum:
\begin{eqnarray*}
\Psi(x)=\sum_{m=1}^{\infty}\frac{(-1)^{m-1}}{m}(e^x-1)^m
\end{eqnarray*}

Notice that for $x=0$ this is identically zero. Taking a derivative:
\begin{eqnarray*}
\Psi^{'}(x)=e^{x}\sum_{m=0}^\infty (-1)^{m}(e^x-1)^{m}\\
=e^{x}e^{-x}=1
\end{eqnarray*}

Consequently, $\Psi(x)=x$, and $a_N$ vanishes for $N \geq 2$, so the first identity we sought is proven. Moving on to the second, define:
\begin{eqnarray*}
b_N=\sum_{m=2}^N \frac{(-1)^{m-1}}{m}\sum_{k_1+...+k_m=N, k_i > 0}{N \choose k_1,...,k_m}
\sum_{i \neq j}^{m} k_ik_j\\
=
\sum_{m=2}^N (-1)^{m-1}(m-1)\sum_{k_1+...+k_m=N, k_i > 0}{N \choose k_1,...,k_m}
k_1k_2
\end{eqnarray*}

Computing the generating function:
\begin{eqnarray*}
\Psi(x)=\sum_{N=1}^{\infty} b_N \frac{x^N}{N!}
=\sum_{N=1}^{\infty} \sum_{m=2}^N (-1)^{m-1}(m-1)\sum_{k_1+...+k_m=N, k_i > 0}{N \choose k_1,...,k_m}
\frac{x^N}{N!} k_1k_2\\
=\sum_{m=2}^\infty (-1)^{m-1}(m-1)\sum_{N=1}^{\infty}\sum_{k_1+...+k_m=N, k_i > 0}{N \choose k_1,...,k_m}
\frac{x^N}{N!} k_1k_2\\
=\sum_{m=2}^\infty (-1)^{m-1}(m-1)\sum_{k_1>0,..., k_m > 0}
\frac{x^{k_1}...x^{k_m}}{(k_1-1)!(k_2-1)!k_3!...k_m!}
\end{eqnarray*}

We substitute in the following identities:
\begin{eqnarray*}
\sum_{k=1}^{\infty} \frac{x^k}{(k-1)!}
=xe^x\\
\sum_{k=1}^{\infty} \frac{x^k}{k!}
=e^x-1
\end{eqnarray*}

This provides:
\begin{eqnarray*}
\Psi(x)=x^2e^{2x}\sum_{m=2}^\infty (-1)^{m-1}(m-1)(e^x-1)^{m-2}\\
=x^2e^{2x}\sum_{m=1}^\infty (-1)^{m}m(e^x-1)^{m-1}
\end{eqnarray*}

Now:
\begin{eqnarray*}
e^x=\frac{d}{dx}(e^x)=\frac{d}{dx} (e^{2x}\sum (-1)^m(e^x-1)^m)\\
2e^{2x}\sum (-1)^m(e^x-1)^m+e^{3x}\sum (-1)^mm(e^x-1)^{m-1}
\end{eqnarray*}

Or else:
\begin{eqnarray*}
-e^{-2x}=\sum (-1)^mm(e^x-1)^{m-1}
\end{eqnarray*}

Which gives $\Psi(x)=x^2$, and consequently $b_N$ vanishes for $N > 2$. This concludes the proof $\P$.

The next lemma we will need is
an integral identity:

\begin{lemma}\label{com1}
Consider the following expression:
\begin{eqnarray*}
\int f(z)e^{-n|z-w|^2}\exp[-2ni(\Re(z)\Im(w)-\Im(z)\Re(w))]\\
\exp[-n(\Re(z)-a)^2/2c]\exp[-n(\Im(z)-b)^2/2d]
\exp[-inp\Re(z)]\exp[-inq\Im(z)]d^2z
\end{eqnarray*}
We may rewrite this integral as
\begin{eqnarray*}
\int
\Psi_{f}(\xi)
\exp\left(-n\left(
\Im(w)-\frac{4c(\Re(\xi)/2n-p/2)+2b}
{4c+2}\right)^2\right)
\exp\left(-n\left(
\Re(w)-\frac{4c(-\Im(\xi)/2n-q/2)+2a}
{4c+2}\right)^2\right)
\\
\exp\left(-in (-2(1/c)b+2p-2\Re(\xi)/n) \left[\frac{\Re(w)}{2+(1/c)}\right]\right)
\exp\left(-in (2(1/c)a+2q-2\Im(\xi)/n) \left[\frac{\Im(w)}{2+(1/c)}\right]\right)d\xi
\end{eqnarray*}

Here we define:
\begin{eqnarray*}
\Psi_{f}(\xi)=\frac{1}{4\pi n}\frac{2c}{1+2c}
\hat{f}(\xi)
\exp\left(-n\frac{((\Re(\xi)/2n-p/2)-b)^2}{\frac{1}{2c}+2+2c}\right)
\exp\left(-n\frac{((-\Im(\xi)/2n+q/2)-b)^2}{\frac{1}{2c}+2+2c}\right)\\
\exp\left(-i (np-\Re(\xi)) \left[\frac{(1/c)a}{2+(1/c)}\right]\right)
\exp\left(-i (nq-\Im(\xi)) \left[\frac{(1/c)b}{2+(1/c)}\right]\right)
\end{eqnarray*}

Here $a,b,c,d$ are constants and $f$ is a Schwartz function.
\end{lemma}

\textbf{Proof:}

First, recall that the product of two Gaussian functions is another Gaussian.

In fact, we have the classical identity:

\begin{eqnarray*}
\frac{1}{\sqrt{2\pi}\sigma_1}\mbox{exp}\left(\frac{(x-\mu_1)^2}{2\sigma_1^2}\right)
\frac{1}{\sqrt{2\pi}\sigma_2}\mbox{exp}\left(\frac{(x-\mu_2)^2}{2\sigma_2^2}\right)\\
=
\frac{1}{\sqrt{2\pi}\sqrt{\sigma_1^2+\sigma_2^2}}\mbox{exp}\left(\frac{(\mu_1-\mu_2)^2}{2(\sigma_1^2+\sigma_2^2)}\right)\\
\times
\frac{1}{\sqrt{2\pi}}
\left(\frac{\sigma_1^2\sigma_2^2}{\sigma_1^2+\sigma_2^2}\right)^{-1/2}
\mbox{exp}\left(\frac{\left(x-\frac{\sigma_1^{-2}\mu_1+\sigma_2^{-2}\mu_2}{\sigma_1^{-2}+\sigma_2^{-2}}\right)^2}{2\frac{\sigma_1^2\sigma_2^2}{\sigma_1^2+\sigma_2^2}}\right)
\end{eqnarray*}

Now consider the expression:
\begin{eqnarray*}
\int f(z)
\mbox{exp}\left(-n|z-w|^2\right)
\mbox{exp}\left(-2ni(\Re(z)\Im(w)-\Im(z)\Re(w))\right)\\
\mbox{exp}\left(-n(\Re(z)-a)^2/2c\right)
\mbox{exp}\left(-n(\Im(z)-b)^2/2d \right)
\mbox{exp}\left(-inp\Re(z)\right)
\mbox{exp}\left(-inq\Im(z)\right)\\
=
\int f(z)\mbox{exp}\left(-n(\Re(z)-\Re(w))^2\right)
\mbox{exp}\left(-n(\Re(z)-a)^2/2c\right)\\
\mbox{exp}\left(-n(\Im(z)-\Im(w))^2\right)
\mbox{exp}\left(-n(\Im(z)-b)^2/2d\right)\\
\mbox{exp}\left(-i\Re(z)[2n\Im(w)+np]\right)
\mbox{exp}\left(-i\Im(z)[-2n\Re(w)+nq]\right)
\end{eqnarray*}

Applying the formula for products of Gaussians:
\begin{eqnarray*}
\mbox{exp}\left(-\frac{(x-\Re(w))^2}{2(1/2n)}\right)
\mbox{exp}\left(-\frac{(x-a)^2}{2c(1/n)}\right)\\
=
\mbox{exp}\left(-\frac{(\Re(w)-a)^2}{2((1/2n)+c(1/n))}\right)
\times
\mbox{exp}\left(-\frac{\left(x-\frac{2n\Re(w)+(n/c)a}{2n+(n/c)}\right)^2}{2\frac{(1/2n)c(1/n)}{(1/2n)+c(1/n)}}\right)\\
=
\mbox{exp}\left(-\frac{n(\Re(w)-a)^2}{1+2c}\right)
\times
\mbox{exp}\left(-n\frac{\left(x-\frac{2\Re(w)+(1/c)a}{2+(1/c)}\right)^2}{\frac{2c}{1+2c}}\right)
\end{eqnarray*}

Substituting this in, our expression is now:
\begin{eqnarray*}
\int f(z)\mbox{exp}\left(-\frac{n(\Re(w)-a)^2}{1+2c}\right)
\times
\mbox{exp}\left(-n\frac{\left(\Re(z)-\frac{2\Re(w)+(1/c)a}{2+(1/c)}\right)^2}{\frac{2c}{1+2c}}\right)
\\
\mbox{exp}\left(-\frac{n(\Im(w)-b)^2}{1+2c}\right)
\times
\mbox{exp}\left(-n\frac{\left(\Im(z)-\frac{2\Im(w)+(1/c)b}{2+(1/c)}\right)^2}{\frac{2c}{1+2c}}\right)\\
\mbox{exp}\left(-i\Re(z)[2n\Im(w)+np]\right)
\times
\mbox{exp}\left(-i\Im(z)[-2n\Re(w)+nq]\right)
\end{eqnarray*}

Taking a Fourier transform, we obtain the equivalent formulation:
\begin{eqnarray*}
\frac{1}{4\pi^2}\frac{\pi}{n}\frac{2c}{1+2c}
\int \hat{f}(\xi)\mbox{exp}\left(-\frac{n(\Re(w)-a)^2}{1+2c}\right)\mbox{exp}\left(-\frac{n(\Im(w)-b)^2}{1+2c}\right)\\
\times
\mbox{exp}\left(-\frac{2c}{1+2c}\frac{\left((2n\Im(w)+np-\Re(\xi))\right)^2}{4n}\right)
\mbox{exp}\left(-i (2n\Im(w)+np-\Re(\xi)) [\frac{2\Re(w)+(1/c)a}{2+(1/c)}]\right)
\\
\times
\mbox{exp}\left(-\frac{2c}{1+2c}\frac{\left((-2n\Re(w)+nq-\Im(\xi))\right)^2}{4n}\right)
\mbox{exp}\left(-i (-2n\Re(w)+nq-\Im(\xi)) [\frac{2\Im(w)+(1/c)b}{2+(1/c)}]\right)
\end{eqnarray*}

Upon rearrangement:
\begin{eqnarray*}
\frac{1}{4\pi n}\frac{2c}{1+2c}
\int \hat{f}(\xi)
\mbox{exp}\left(-\frac{2nc}{1+2c}\left((\Im(w)+p/2-\Re(\xi)/2n)\right)^2\right)
\times
\mbox{exp}\left(-\frac{n(\Im(w)-b)^2}{1+2d}\right)
\\
\mbox{exp}\left(-\frac{2nc}{1+2c}\left((\Re(w)-q/2+\Im(\xi)/2n)\right)^2\right)
\times
\mbox{exp}\left(-\frac{n(\Re(w)-a)^2}{1+2c}\right)
\\
\mbox{exp}\left(-i (2n\Im(w)+np-\Re(\xi)) [\frac{2\Re(w)+(1/c)a}{2+(1/c)}]\right)\\
\mbox{exp}\left(-i (-2n\Re(w)+nq-\Im(\xi)) [\frac{2\Im(w)+(1/d)b}{2+(1/d)}]\right)
\end{eqnarray*}

The formula for products of Gaussians implies:
\begin{eqnarray*}
\mbox{exp}\left(-\frac{(\Im(w)-(\Re(\xi)/2n-p/2))^2}{2\frac{1+2c}{4nc}}\right)
\mbox{exp}\left(-\frac{(\Im(w)-b)^2}{2\frac{1+2d}{2n}}\right)\\
=
\mbox{exp}\left(-\frac{((\Re(\xi)/2n-p/2)-b)^2}{2(\frac{1+2c}{4nc}+\frac{1+2d}{2n})}\right)
\times
\mbox{exp}\left(-\frac{\left(\Im(w)-\frac{\frac{4nc}{1+2c}(\Re(\xi)/2n-p/2)+\frac{2n}{1+2d}b}{\frac{4nc}{1+2c}+
\frac{2n}{1+2d}}\right)^2}{2\frac{\frac{1+2c}{4nc}\frac{1+2d}{2n}}{\frac{1+2c}{4nc}+\frac{1+2d}{2n}}}\right)
\\
=
\mbox{exp}\left(-n\frac{((\Re(\xi)/2n-p/2)-b)^2}{\frac{1}{2c}+2+2c}\right)
\times
\mbox{exp}\left(-n\left(
\Im(w)-\frac{4c(\Re(\xi)/2n-p/2)+2b}
{4c+2}\right)^2\right)
\end{eqnarray*}

Similarly, we may also rewrite:
\begin{eqnarray*}
\mbox{exp}\left(-i (2n\Im(w)+np-\Re(\xi)) [\frac{2\Re(w)+(1/c)a}{2+(1/c)}]\right)\\
\mbox{exp}\left(-i (-2n\Re(w)+nq-\Im(\xi)) [\frac{2\Im(w)+(1/d)b}{2+(1/c)}]\right)\\
=
\mbox{exp}\left(-i (-n(1/c)b+np-\Re(\xi)) [\frac{2\Re(w)}{2+(1/c)}]\right)\\
\mbox{exp}\left(-i (np-\Re(\xi)) [\frac{(1/c)a}{2+(1/c)}]\right)\\
\mbox{exp}\left(-i (n(1/c)a+nq-\Im(\xi)) [\frac{2\Im(w)}{2+(1/c)}]\right)\\
\mbox{exp}\left(-i (nq-\Im(\xi)) [\frac{(1/c)b}{2+(1/c)}]\right)
\end{eqnarray*}

Plugging both of these equalities in, our original expression becomes:
\begin{eqnarray*}
\frac{1}{4\pi n}\frac{2c}{1+2c}
\int \hat{f}(\xi)
\mbox{exp}\left(-n\frac{((\Re(\xi)/2n-p/2)-b)^2}{\frac{1}{2c}+2+2c}\right)
\times
\mbox{exp}\left(-n\left(
\Im(w)-\frac{4c(\Re(\xi)/2n-p/2)+2b}
{4c+2}\right)^2\right)
\\
\mbox{exp}\left(-n\frac{((-\Im(\xi)/2n+q/2)-b)^2}{\frac{1}{2c}+2+2c}\right)
\times
\mbox{exp}\left(-n\left(
\Re(w)-\frac{4c(-\Im(\xi)/2n-q/2)+2a}
{4c+2}\right)^2\right)
\\
\mbox{exp}\left(-in (-2(1/c)b+2p-2\Re(\xi)/n) [\frac{\Re(w)}{2+(1/c)}]\right)\\
\mbox{exp}\left(-in (2(1/c)a+2q-2\Im(\xi)/n) [\frac{\Im(w)}{2+(1/c)}]\right)\\
\mbox{exp}\left(-i (np-\Re(\xi)) [\frac{(1/c)a}{2+(1/c)}]\right)
\mbox{exp}\left(-i (nq-\Im(\xi)) [\frac{(1/c)b}{2+(1/c)}]\right)
\end{eqnarray*}

This is as we sought.
$\P$

Applying this identity to nested integrals, we obtain:

\begin{lemma}\label{com2}
The following integral identity holds:
\begin{eqnarray*}
\int f_1(z_1)....f_m(z_m)\exp(-n|z_1-z_2|^2)
...\exp(-n|z_m-z_1|^2)\\
\exp(2ni[\Re(z_1)\Im(z_2)-\Im(z_1)\Re(z_2)])
...
\exp(2ni[\Re(z_n)\Im(z_1)-\Im(z_n)\Re(z_1)])
d^2(z_1)...d^2(z_m)
\\
=
\left(\frac{\pi}{2n}\right)^{m-1}\int \left(\prod_{j=1}^m \hat{f}(z)\right)d^2z
-\frac{1}{2n}\left(\frac{\pi}{2n}\right)^{m-1}\sum_{j < k} \int \frac{df_j}{dz}\frac{df_k}{dz}\prod_{h\neq j,k} f(z)d^2z
+O(1/n)
\end{eqnarray*}
Here $f_i$ is a family of Schwartz functions.
\end{lemma}

\textbf{Proof:}

We apply lemma \ref{com1} to the inner most integral, setting:

\begin{eqnarray*}
a=\Re(z); b=\Im(z); p=-2\Im(z)\\ q=2\Re(z); c=1/2
\end{eqnarray*}

Moving on to the next integral, we apply lemma \ref{com1} again, this time with:
\begin{eqnarray*}
a= -\Im(\xi_1)/4n+\Re(z)\\
b = \Re(\xi_1)/4n+\Im(z)\\
q=2a, p=-2b\\
c=1/2
\end{eqnarray*}

The pattern continues, and after evaluating all integrals but the last in this fashion we obtain:
\begin{eqnarray*}
\left(\frac{1}{8\pi n}\right)^{m-1}\int f_m(z)\left(\prod_{j=1}^{m-1} \hat{f_j}(\xi_j) \right)\\
\mbox{exp}\left(-n(\Re(z)-\Re(z)-\sum_{j=1}^m \Im(\xi_j)/4n)^2\right)
\mbox{exp}\left(-n(\Im(z)-\Im(z)+\sum_{j=1}^m \Re(\xi_j)/4n)^2\right)\\
\mbox{exp}\left(-i2n(\Re(z)-\sum_{j=1}^m \Im(\xi_j)/4n)\Im(z)\right)
\mbox{exp}\left( i2n(\Im(z)+\sum_{j=1}^m \Re(\xi_j)/4n)\Re(z)\right)\\
\times
\mbox{exp}\left( -\frac{\sum_{j=1}^m(\Re(\xi_j))^2}{16n}\right)
\mbox{exp}\left( -\frac{\sum_{j=1}^m(\Im(\xi_j))^2}{16n}\right)\\
\mbox{exp}\left(\frac{i}{2}
\Re(z)[\sum_{j=1}^m \Re(\xi_j)]-\sum_{k=1}^m(\Im(\xi_k)/4n)[\sum_{j>k}^{m} \Re(\xi_j)]
\right)\\
\mbox{exp}\left(\frac{i}{2}
\Im(z)[\sum_{j=1}^m \Im(\xi_j)]+\sum_{k=1}^m(\Re(\xi_k)/4n)[\sum_{j>k}^{m} \Im(\xi_j)]
\right)d^2z d^2\xi_1...d^2\xi_{m-1}
\end{eqnarray*}

Or:
\begin{eqnarray*}
\left(\frac{1}{8\pi n}\right)^{m-1}\int f_m(z)\left(\prod_{j=1}^{m-1} \hat{f_j}(\xi_j) \right)\\
\mbox{exp}\left( -\frac{[\sum_{j=1}^m(\Re(\xi_j)]^2+\sum_{j=1}^m(\Re(\xi_j))^2}{16n}\right)
\mbox{exp}\left( -\frac{[\sum_{j=1}^m(\Im(\xi_j)]^2+\sum_{j=1}^m(\Im(\xi_j))^2}{16n}\right)\\
\mbox{exp}\left(-\frac{i}{8n}\sum_{k=1}^m \Im(\xi_k)\sum_{j>k}^{m} \Re(\xi_j)
\right)
\mbox{exp}\left(\frac{i}{8n}
\sum_{k=1}^m \Re(\xi_k)\sum_{j>k}^{m} \Im(\xi_j)
\right)
\\
\mbox{exp}\left(i
\Re(z) \sum_{j=1}^m \Re(\xi_j)
\right)
\mbox{exp}\left(i
\Im(z) \sum_{j=1}^m \Im(\xi_j)
\right)d^2z d^2\xi_1...d^2\xi_{m-1}
\end{eqnarray*}

Or, after further rearragement:
\begin{eqnarray*}
\left(\frac{1}{8\pi n}\right)^{m-1}\int f_m(z)\left(\prod_{j=1}^m \hat{f_j}(\xi_j) \right)
\mbox{exp}\left(i
\Re(z) \sum_{j=1}^m \Re(\xi_j)
\right)
\mbox{exp}\left(i
\Im(z) \sum_{j=1}^m \Im(\xi_j)
\right)\\
\mbox{exp}\left( -\frac{\sum_{j=1}^m \Re(\xi_j)^2+\sum_{j\neq k}^m \Re(\xi_j)\Re(\xi_k)}{8n}\right)
\mbox{exp}\left( -\frac{\sum_{j=1}^m \Im(\xi_j)^2+\sum_{j\neq k}^m \Im(\xi_j)\Im(\xi_k)}{8n}\right)
\\
\mbox{exp}\left(\frac{i}{8n}\Big(-\sum_{k=1}^m \Im(\xi_k)[\sum_{j>k}^m\Re(\xi_k)]+\sum_{k=1}^m\Re(\xi_k)[\sum_{j>k}\Im(\xi_j)]
\Big)
\right)
d^2z d^2\xi_1...d^2\xi_{m-1}
\end{eqnarray*}

Taking a Fourier transform, our expression becomes the more tractable:

\begin{eqnarray*}
\left(\frac{1}{8\pi n}\right)^{m-1}\int \left(\prod_{j=1}^m \hat{f}(\xi_j) \right) \overline{\hat{f}\left(\sum_{j=1}^m \xi_j\right)}
\mbox{exp}\left(-\frac{1}{8n}
\sum_{j=1}^m |\xi_j|^2
\right)
\mbox{exp}\left(-\frac{1}{8n}
\sum_{j> k} \xi_{k}\overline{\xi_{j}}
\right)
\end{eqnarray*}

We have:
\begin{eqnarray*}
\int \hat{f}_1(z_1)...\hat{f}_m(z_{m-1})\overline{\hat{f}_{m}(z_1+...+z_{m-1})}|z_1|^2\\
=\int \hat{f}_1(z_1)\overline{z_1}...\hat{f}_m(z_{m-1})\overline{\hat{f}_{m}(z_1+...+z_{m-1})\overline{\sum_{j=1} z_j}}\\
-\sum_{j\neq 1}\hat{f}_1(z_1)\overline{z_1}\hat{f}(z_j)z_j...\hat{f}_m(z_{m-1})\overline{\hat{f}_{m}(z_1+...+z_{m-1})}
\end{eqnarray*}

Consequently:
\begin{eqnarray*}
\sum_j \int \hat{f}_1(z_1)...\hat{f}_m(z_{m-1})\overline{\hat{f}_{m}(z_1+...+z_{m-1})}|z_j|^2\\
=\sum_j \int \hat{f}_j(z_j)\overline{z_j}\left(\prod_{k\neq j}\hat{f}_k(z_k)\right)\overline{\hat{f}_{m}(z_1+...+z_{m-1})\overline{\sum_{j=1} z_j}}\\
-\sum_{j \neq k}\int \hat{f}_j(z_j)\overline{z_j} \hat{f}_k(z_k)z_k \left(\prod_{h \neq j,k} \hat{f}_h(z_h)
\right)\overline{\hat{f}_{m}(z_1+...+z_{m-1})}
\end{eqnarray*}

Then:
\begin{eqnarray*}
\sum_j \int \hat{f}_1(z_1)...\hat{f}_m(z_{m-1})\overline{\hat{f}_{m}(z_1+...+z_{m-1})}|z_j|^2\\
+\sum_{j <k}\int \hat{f}_1(z_1)...\hat{f}_m(z_{m-1})\overline{\hat{f}_{m}(z_1+...+z_{m-1})}z_j\overline{z_k}\\
=
\sum_j \int \left[\hat{f}_j(z_j)i\overline{z_j}\right]\left(\prod_{k\neq j}\hat{f}_k(z_k)\right)\overline{\hat{f}_{m}(z_1+...+z_{m-1})i\overline{\sum_{j=1} z_j}}\\
+\sum_{j < k}\int \left[ \hat{f}_j(z_j)i\overline{z_j}\right] \left[\hat{f}_k(z_k)iz_k\right] \left(\prod_{h \neq j,k} \hat{f}_h(z_h)\right)
\overline{\hat{f}_{m}(z_1+...+z_{m-1})}
\end{eqnarray*}

By Plancharel's theorem, this is:
\begin{eqnarray*}
4(4\pi^2)^{m-1}\sum_{j < k} \int \frac{df_j}{dz}\frac{df_k}{dz}\prod_{h\neq j,k} f(z)dz
\end{eqnarray*}

By a Taylor expansion then:

\begin{eqnarray*}
\left(\frac{1}{8\pi n}\right)^{m-1}\int \left(\prod_{j=1}^m \hat{f}(\xi_j) \right) \overline{\hat{f}\left(\sum_{j=1}^m \xi_j\right)}
\mbox{exp}\left(-\frac{1}{8n}
\sum_{j=1}^m |\xi_j|^2
\right)
\mbox{exp}\left(-\frac{1}{8n}
\sum_{k \neq j} \xi_{k}\overline{\xi_{j}}
\right)\\
=\left(\frac{\pi}{2n}\right)^{m-1}\int \left(\prod_{j=1}^m \hat{f}(z)\right)d^2z
-\frac{1}{2n}\left(\frac{\pi}{2n}\right)^{m-1}\sum_{j < k} \int \frac{df_j}{dz}\frac{df_k}{dz}\prod_{h\neq j,k} f(z)d^2z
+O(1/n)
\end{eqnarray*}

Substituting these formulas in, the lemma is proven. $\P$

With these lemmas in hand, we may now prove the main result of this section.

\textbf{Proof of Proposition \ref{highc}:}

Recall our expression for the $k$-th cumulant:

\begin{eqnarray*}
\mathcal{C}_k(\Lambda_{2n}[f])=\sum_{m=1}^k\frac{(-1)^{m-1}}{m}\sum_{k_1+...+k_m=k}\frac{k!}{k_1!...k_m!}\\
\times
\int_{C^k}\left(\prod_{h=1}^m f(z_h)^{k_h}\right)S_{2n}(z_1,z_2)...S_{2n}(z_m,z_1)d^2(z_1)...d^2(z_m)
\end{eqnarray*}

To understand this expression, we will first concentrate on simplifying the integral. To this end, fix a Schwartz function $\psi$, and assume $\sum_{j=1}^m q_j=M$, for some positive integer $M$ and collection of integers $q_j > 0$. (For convenience, we also set $q_{m+1}=q_1$.)

We will begin by investigating the following expression:
\begin{eqnarray*}
(2n)^m\int \psi^{q_1}(z_1)...\psi^{q_m}(z_m)S_n(\sqrt{2n}z_1,\sqrt{2n}z_2)...S_n(\sqrt{2n}z_m,\sqrt{2n}z_1)\\
=\left(\frac{-i(2n)}{\sqrt{2\pi}}\right)^m\left(\sqrt{2n}\right)^m \int \psi^{q_1}(x_1)...\psi^{q_m}(x_m)
e^{-n\sum_{j}(z_j-\overline{z_{j+1}})^2}\prod_j \left(z_j-\overline{z_{j+1}}\right)
\prod_j\left(\mbox{erfc}(2\sqrt{n}\Im(z_j)\right)
\end{eqnarray*}

Applying the asymptotic expansion of the complementary error function, this expression becomes:
\begin{eqnarray*}
\left(\frac{-i(2n)}{\sqrt{2\pi}}\right)^m
\left(\frac{\sqrt{2n}}{2\sqrt{\pi n}}\right)^m
\int \left(\prod_j\psi^{q_j}(x_j)\left(
1-\frac{1}{8n\Im(z_j)^2}+O(n^{-2})
\right)\right)
\\
\prod_j \left(\frac{1}{\Im(z_j)}(z_j-\overline{z_{j+1}})\right)
\mbox{exp}\left(
-n[|z_1-z_2|^2+...+|z_{m-1}-z_{m}|^2+|z_m-z_1|^2]
\right)\\
\mbox{exp}\left(
-2ni[(\Re(z_1)\Im(z_2)-\Re(z_2)\Im(z_1))+...+
(\Re(z_m)\Im(z_1)-\Re(z_1)\Im(z_m))]
\right)
\end{eqnarray*}
 
In order to pick out first order (and second order) terms featured in this expression, we need to expand the following product:
\begin{eqnarray*}
\prod_j \left(\frac{1}{\Im(z_j)}(z_j-\overline{z_{j+1}})\right)
=
\prod_j \left(\frac{1}{\Im(z_j)}([z_j-z_{j+1}]+2i\Im(z_{j+1}))\right)
\end{eqnarray*}

At leading order, we have the term corresponding to exclusively factors of the form $2i\Im(z_{j+1})$. Define:
\begin{eqnarray*}
T_1=\left(\frac{2n}{\pi}\right)^m
\int \left(\prod_j\psi^{q_j}(x_j)\left(
1-\frac{1}{8n\Im(z_j)^2}+O(n^{-2})
\right)\right)
\\
\mbox{exp}\left(
-n[|z_1-z_2|^2+...+|z_{m-1}-z_{m}|^2+|z_m-z_1|^2]
\right)\\
\mbox{exp}\left(
-2ni[(\Re(z_1)\Im(z_2)-\Re(z_2)\Im(z_1))+...+
(\Re(z_m)\Im(z_1)-\Re(z_1)\Im(z_m))]
\right)
\end{eqnarray*}

We will see that $T_1$  is of order
$n$. We will also need to consider second order terms (in which a single factor in the above product is not of the form  $2i\Im(z_{j+1})$). This motivates the next definition:

\begin{eqnarray*}
T_2=\sum_k
\left(\frac{-ni}{\pi}\right)^m
\int \left(\prod_j\psi^{q_j}(x_j)\right)
\frac{(2i)^{m-1}}{\Im(z_k)}(z_{k+1}-z_k)
\\
\mbox{exp}\left(
-n[|z_1-z_2|^2+...+|z_{m-1}-z_{m}|^2+|z_m-z_1|^2]
\right)\\
\mbox{exp}\left(
-2ni[(\Re(z_1)\Im(z_2)-\Re(z_2)\Im(z_1))+...+
(\Re(z_m)\Im(z_1)-\Re(z_1)\Im(z_m))]
\right)
\end{eqnarray*}

We will see that $T_2$ is of order unity. All other, smaller terms vanish in the limit as $n$ becomes large.

Now, $T_1$ is easily handled by Lemma \ref{com2}. To this end, set:

\begin{eqnarray*}
f_{j}(x)=\psi(x_j)^{q_j}\left(
1-\frac{1}{8n \Im(x_j)^2}
\right)
\end{eqnarray*}

The lemma now yields that, up to small terms which vanish in the limit,  we have:

\begin{eqnarray*}
T_1\approx
\left(\frac{\pi}{2n}\right)^{m-1}\left(\frac{2n}{\pi}\right)^m\int \psi(z)^M d^2z -
\left(\frac{\pi}{2n}\right)^{m-1}\frac{m2^{m-3}n^{m-1}}{\pi^m}\int \frac{\psi(z)^M}{\Im(z)^2} d^2z \\
+\left(\frac{\pi}{2n}\right)^{m-1}\left(\frac{2n}{\pi}\right)^m\frac{1}{8n}\sum_{j < k}q_jq_k
\int |\psi^{'}(z)|^2\psi^{M-2}(z)d^2z
\\
\approx
\frac{2n}{\pi}\int \psi(z)^M d^2z -
\frac{m}{4\pi}\int \frac{\psi(z)^M}{\Im(z)^2} d^2z \\
+\frac{1}{4\pi}\sum_{j < k}q_jq_k
\int |\psi^{'}(z)|^2\psi^{M-2}(z)d^2z
\end{eqnarray*}

Now, we turn to evaluate $T_2$. By Lemma \ref{com2}:

\begin{eqnarray*}
\left(\frac{-ni}{\pi}\right)^m
\int \left(\prod_j\psi^{q_j}(x_j)\right)
\frac{(2i)^{m-1}}{\Im(z_1)}(z_{1})
\\
\mbox{exp}\left(
-n[|z_1-z_2|^2+...+|z_{m-1}-z_{m}|^2+|z_m-z_1|^2]
\right)\\
\mbox{exp}\left(
-2ni[(\Re(z_1)\Im(z_2)-\Re(z_2)\Im(z_1))+...+
(\Re(z_m)\Im(z_1)-\Re(z_1)\Im(z_m))]
\right)\\
=
\left(\frac{\pi}{2n}\right)^{m-1}\left(\frac{-ni}{\pi}\right)^m
(2i)^{m-1}\int \frac{z\psi^{M}(z)}{\Im(z)}dz
-\left(\frac{\pi}{2n}\right)^{m-1}\frac{1}{2n}\left(\frac{-ni}{\pi}\right)^m
(2i)^{m-1}\\
\times
\Big[\sum_{k_1<k_2}\int q_{k_1}q_{k_2}\frac{z\left(\frac{d\psi}{dz}\right)^2\psi^{M-2}(z)}{\Im(z)}dz
+\sum_{k \neq 1}\int q_{k}\left(
\frac{2}{\Im(z)}+\frac{iz}{\Im(z)^2}
\right)
\left(\frac{d\psi}{dz}\right)\psi^{M-1}(z)dz
\end{eqnarray*}

And similarly:

\begin{eqnarray*}
\left(\frac{-ni}{\pi}\right)^m
\int \left(\prod_j\psi^{q_j}(x_j)\right)
\frac{(2i)^{m-1}}{\Im(z_1)}(z_{2})
\\
\mbox{exp}\left(
-n[|z_1-z_2|^2+...+|z_{m-1}-z_{m}|^2+|z_m-z_1|^2]
\right)\\
\mbox{exp}\left(
-2ni[(\Re(z_1)\Im(z_2)-\Re(z_2)\Im(z_1))+...+
(\Re(z_m)\Im(z_1)-\Re(z_1)\Im(z_m))]
\right)\\
=
\left(\frac{\pi}{2n}\right)^{m-1}\left(\frac{-ni}{\pi}\right)^m
(2i)^{m-1}\int \frac{z\psi^{M}(z)}{\Im(z)}dz
-\left(\frac{\pi}{2n}\right)^{m-1}\frac{1}{4n}\left(\frac{-ni}{\pi}\right)^m
(2i)^{m-1}\\
\times
\Big[\sum_{k_1<k_2}\int q_{k_1}q_{k_2}\frac{z\left(\frac{d\psi}{dz}\right)^2\psi^{M-2}(z)}{\Im(z)}dz
+2\sum_{k \neq 2}\int q_{k}\frac{\left(\frac{d\psi}{dz}\right)\psi^{M-1}(z)}{\Im(z)}dz\\
+i\sum_{k \neq 1}\int q_k \frac{z\left(\frac{d\psi}{dz}\right)\psi^{M-1}(z)}{\Im(z)^2}dz
-2i\int \frac{\psi^{M}(z)}{\Im(z)^2}dz\Big]
\end{eqnarray*}

Consequently:
\begin{eqnarray*}
\left(\frac{-ni}{\pi}\right)^m
\int \left(\prod_j\psi^{q_j}(x_j)\right)
\frac{(2i)^{m-1}}{\Im(z_1)}(z_{2}-z_1)
\\
\mbox{exp}\left(
-n[|z_1-z_2|^2+...+|z_{m-1}-z_{m}|^2+|z_m-z_1|^2]
\right)\\
\mbox{exp}\left(
-2ni[(\Re(z_1)\Im(z_2)-\Re(z_2)\Im(z_1))+...+
(\Re(z_m)\Im(z_1)-\Re(z_1)\Im(z_m))]
\right)\\
=
\frac{1}{4}
\frac{1}{\pi}(-1)^m i^{2m-1}
\times
\Big[i\int \frac{\psi^{M}(z)}{\Im(z)^2}dz
\\
-2
\sum_{k \neq 2}\int q_{k}\frac{\psi^{'}(z)\psi^{M-1}(z)}{\Im(z)}dz
+2\sum_{k \neq 1}\int q_{k}
\frac{\psi^{'}(z)\psi^{M-1}(z)}{\Im(z)}dz
\Big]
\end{eqnarray*}

Applying the analagous formula for every summand in $T_2$:

\begin{eqnarray*}
T_2=\frac{m}{4\pi}
\int \frac{\psi^{M}(z)}{\Im(z)^2}dz
\end{eqnarray*}

This provides, up to small terms which vanish in the limit, the following expansion:
\begin{eqnarray*}
T_1+T_2\approx
\frac{2n}{\pi}\int \psi(z)^M d^2z
+\frac{1}{4\pi}\sum_{j < k}q_jq_k
\int |\psi^{'}(z)|^2\psi^{M-2}(z)d^2z
\end{eqnarray*}

Comparing with our Costin-Lebowtiz type formula for the cumulants of the linear statistic, we obtain:
\begin{eqnarray*}
\mathcal{C}_k(\Lambda_{2n}[f])=\sum_{m=1}^k\frac{(-1)^{m-1}}{m}\sum_{k_1+...+k_m=k}\frac{k!}{k_1!...k_m!}\\
\Big(\frac{2n}{\pi}\int \psi(z)^M d^2z
+\frac{1}{4\pi}\sum_{j < i}k_j k_i
\int |\psi^{'}(z)|^2\psi^{M-2}(z)d^2z+O(n^{-\epsilon})\Big)
\end{eqnarray*}

Applying Lemma \ref{com0}, we see that the contributions of the $O(n)$ and $O(1)$ terms are both zero. Consequently, $\mathcal{C}_k(\Lambda_{2n}[f])$
vanishes in the limit, and the proposition is proven $\P$.

\subsection{Variance in the Real/Real Case}

Now we consider the case of test functions supported on the real line itself. First, we will compute the limiting variance:

\begin{proposition}\label{VarR}
Let $f$ be a real valued test function, compactly supported on $(0,1)$, and let $\lambda_j$ denote eigenvalues drawn from the real Ginibre ensemble.
Then normalized variance of the random variable $S_n(f)$ satisfies the following identity:
\begin{eqnarray*}
\lim_{n\to \infty} \text{\emph{Var}}\left[\frac{1}{n^{1/4}}S_n[f]\right]=\left(\frac{2-\sqrt{2}}{\sqrt{\pi}}\right)\int f(x)^2 dx
\end{eqnarray*}
In fact, it is sufficient to assume that $f$ has one continuous derivative.
\end{proposition}

Previously, we have computed the following formula for the variance:

\begin{eqnarray*}
\mbox{Var}(S_{2n}[f])=
\sqrt{2n}\int f(x)^2S_{2n}(\sqrt{2n}x,\sqrt{2n}x)dx-
2n\int f(x)f(y)\\
\times\left(
D_{2n}(\sqrt{2n}x,\sqrt{2n}y)I_{2n}(\sqrt{2n}x,\sqrt{2n}y)
+S_{2n}(\sqrt{2n}x,\sqrt{2n}y)S_{2n}(\sqrt{2n}y,\sqrt{2n}x)
\right)dxdy
\end{eqnarray*}

To evaluate this expression, we will make repeated use of the following result of Borodin and Sinclair:
\begin{lemma}
\label{approxim}
For real valued $u$ with $|u|<1$, we have:
\begin{eqnarray*}\label{approx}
e^{-2nu}\sum_{m=0}^{2n-2}\frac{(2nu)^m}{m!}=1+O(1/n)
\end{eqnarray*}
The implied constant depends on the distance from $u$ to the edge, $\pm 1$.
\end{lemma}

We will also need the following limit, which follows from the proof of Lemma 9.3 in Borodin and Sinclair.

\begin{lemma}\label{rsmall}
For real valued $x$ and $y$ with modulus less than one, we have:
\begin{eqnarray*}
\lim_{n\to \infty} \frac{1}{\sqrt{2\pi}}\frac{1}{(2n-2)!}e^{-y^2/2}2^{n-3/2}\emph{sgn}(y)
y^{2n-1}\gamma(n-\frac{1}{2},\frac{x^2}{2})=0
\end{eqnarray*}
In fact, this expression vanishes exponentially fast, with constants depending on the distance between $\max\{|x|,|y|\}$ and 1.
\end{lemma}

This lemma is chiefly interesting to us in that it provides the following estimate:
\begin{eqnarray*}
S_{2n}(x,y)\approx \frac{1}{\sqrt{2\pi}}e^{-(1/2)(x-y)^2}e^{-xy}\sum_{m=0}^{2n-2}\frac{(xy)^m}{m!}
\end{eqnarray*}

We are now in a position to prove Proposition \ref{VarR}. As before, we will compute the asymptotic behavior of the
 various terms which comprise our expression for the variance individually, with the results recorded by a sequence at lemmas.

We begin with the following easy lemma:

\begin{lemma}
In the real/real case, we have the following pair of asymptotic formulas:
\begin{eqnarray*}
\sqrt{2n}\int f(x)^2S_{2n}(\sqrt{2n}x,\sqrt{2n}x)dx=\frac{\sqrt{n}}{\sqrt{\pi}}\int f(x)^2 dx +O\left(\frac{1}{\sqrt{n}}\right)
\end{eqnarray*}

And:
\begin{eqnarray*}
2n\int f(x)f(y)S_{2n}(\sqrt{2n}x,\sqrt{2n}y)S_{2n}(\sqrt{2n}y,\sqrt{2n}x)dxdy\\=
\frac{\sqrt{n}}{\sqrt{2\pi}} \int f(x)^2 dx+O(1/\sqrt{n})
\end{eqnarray*}

The implied constants are allowed to depend on the distance between the support of $f$ and the edge (the points 1 and -1).
\end{lemma}

\textbf{Proof:}

Invoking Lemma \ref{approxim} and Lemma \ref{rsmall}, the first claim reduces to a trivial calculation:

\begin{eqnarray*}
\sqrt{2n}\int f(x)^2S_{2n}(\sqrt{2n}x,\sqrt{2n}x)dx\\
=\frac{\sqrt{n}}{\sqrt{\pi}}\int f(x)^2\left(1+O\left(\frac{1}{n}\right)\right)dx\\
=\frac{\sqrt{n}}{\sqrt{\pi}}\int f(x)^2 dx +||f||^2_{\infty}O\left(\frac{1}{\sqrt{n}}\right)
\end{eqnarray*}

We turn to the second claimed formula. By Lemma \ref{approxim}:
\begin{eqnarray*}
2n\int f(x)f(y)S_{2n}(\sqrt{2n}x,\sqrt{2n}y)S_{2n}(\sqrt{2n}y,\sqrt{2n}x)dxdy\\
\approx 2n\int f(x)f(y)\frac{1}{2\pi}e^{-2n(x-y)^2}dxdy
\end{eqnarray*}

By Plancharel's theorem, this quantity can be rewritten:
\begin{eqnarray*}
\frac{1}{2\pi}\frac{n}{\pi}\int\overline{F[f(x)](\xi)}F[f(x)*e^{-2nx^2}](\xi)d\xi\\
=\frac{1}{2\pi}\frac{n}{\pi}\int \overline{\hat{f}(\xi)}\hat{f}(\xi)\sqrt{\frac{\pi}{2n}}e^{-\xi^2/8n}\\
=\frac{1}{2\pi}\frac{\sqrt{n}}{\sqrt{2\pi}}\int \overline{\hat{f}(\xi)}\hat{f}(\xi)e^{-\xi^2/8n}d\xi
\end{eqnarray*}

By a Taylor expansion, at leading order we have:
\begin{eqnarray*}
\frac{1}{2\pi}\frac{\sqrt{n}}{\sqrt{2\pi}} \int \overline{\hat{f}(\xi)}\hat{f}(\xi) d\xi
=\frac{\sqrt{n}}{\sqrt{2\pi}} \int f(x)^2 dx
\end{eqnarray*}

This concludes the proof. $\P$

Lastly, we show:
\begin{lemma}
In the real/real case, we have the following limit:
\begin{eqnarray*}
2n\int f(x)f(y) D_{2n}(\sqrt{2n}x,\sqrt{2n}y)I_{2n}(\sqrt{2n}x,\sqrt{2n}y)dxdy\\
=
\left(\frac{\sqrt{n}}{\sqrt{2\pi}}-\frac{\sqrt{n}}{\sqrt{\pi}}\right)\int f(x)^2 dx+O\left(\frac{1}{n^{1/4}}\right)
\end{eqnarray*}
Here, $f$ is a smooth test function compactly supported in the open interval $(-1,1)$.
\end{lemma}

\textbf{Proof:}

Recall the following definitions:

\begin{eqnarray*}
D_{2n}(x,y)=(y-x)\frac{1}{\sqrt{2\pi}}e^{-\frac{1}{2}(x-y)^2}e^{-xy}\sum_{m=0}^{2n-2}\frac{(xy)^m}{m!}
\\
I_{2n}(x,y)=e^{-x^2/2}\frac{1}{2\sqrt{\pi}}\int_0^{y^2/2}\frac{e^{-t}}{\sqrt{t}}
\sum_{m=0}^{n-1}\frac{(\sqrt{2t}x)^{2m}}{(2m)!}dt\\-
e^{-y^2/2}\frac{1}{2\sqrt{\pi}}\int_0^{x^2/2}\frac{e^{-t}}{\sqrt{t}}
\sum_{m=0}^{n-1}\frac{(y\sqrt{2t})^{2m}}{(2m)!} dt
+\frac{1}{2}\mbox{sgn}(x-y)
\end{eqnarray*}

Substituting these in:
\begin{eqnarray*}
2n\int f(x)f(y) D_{2n}(\sqrt{2n}x,\sqrt{2n}y)I_{2n}(\sqrt{2n}x,\sqrt{2n}y)dxdy\\
\approx
\frac{2\sqrt{n}n\sqrt{2}}{\sqrt{2\pi}}\int f(x)f(y)(y-x)e^{-n(x-y)^2}\\
\Big[e^{-nx^2}\frac{1}{2\sqrt{\pi}}\int_0^{ny^2}\frac{e^{-t}}{\sqrt{t}}
\sum_{m=0}^{n-1}\frac{(2x\sqrt{tn})^{2m}}{(2m)!}dt\\-
e^{-ny^2}\frac{1}{2\sqrt{\pi}}\int_0^{nx^2}\frac{e^{-t}}{\sqrt{t}}
\sum_{m=0}^{n-1}\frac{(2y\sqrt{nt})^{2m}}{(2m)!} dt\Big]dxdy\\
-
\frac{n\sqrt{n}}{\sqrt{\pi}}\int f(x)f(y)|y-x|e^{-n(x-y)^2}
\end{eqnarray*}

By symmetry, this is:

\begin{eqnarray*}
\frac{4n\sqrt{n}}{\sqrt{\pi}}\int f(x)f(y)(y-x)e^{-n(x-y)^2}\\
\left[e^{-nx^2}\frac{1}{2\sqrt{\pi}}\int_0^{ny^2}\frac{e^{-t}}{\sqrt{t}}
\sum_{m=0}^{n-1}\frac{(2\sqrt{n}x\sqrt{t})^{2m}}{(2m)!}dt\right]dxdy\\
-
\frac{n\sqrt{n}}{\sqrt{\pi}}\int f(x)f(y)|y-x|e^{-n(x-y)^2}
\end{eqnarray*}

The last term is not difficult to deal with. We have the trivial identity:
\begin{eqnarray*}
n\int |x|e^{-nx^2}dx=1
\end{eqnarray*}

Consequently, by a standard approximation to the identity argument, one would expect an estimate such as the following to hold:
\begin{eqnarray*}
\sqrt{n}\frac{n}{\sqrt{\pi}}\int f(x)f(y)|y-x|e^{-n(x-y)^2}dydx
\approx\frac{\sqrt{n}}{\sqrt{\pi}}\int f(x)^2 dx
\end{eqnarray*}

Indeed, this is exactly what happens. We can write:
\begin{eqnarray*}
\frac{n^{3/2}}{\sqrt{\pi}}\int f(y)|y-x|e^{-n(x-y)^2}dy=\frac{n^{3/2}}{\sqrt{\pi}}\int_0^\infty f(y+x) ye^{-ny^2}dy
-\frac{n^{3/2}}{\sqrt{\pi}}\int_{-\infty}^0  f(y+x) ye^{-ny^2}dy\\
=\frac{\sqrt{n}}{\sqrt{\pi}}f(x)+\frac{1}{2}\sqrt{n}\int_0^\infty f^{\prime}(y+x) e^{-ny^2}dy
-\frac{1}{2}\sqrt{n}\int_{-\infty}^0  f^{\prime}(y+x) e^{-ny^2}dy
\end{eqnarray*}

We may also expand:
\begin{eqnarray*}
\frac{2\sqrt{n}}{\sqrt{\pi}}\int_0^\infty f^{\prime}(y+x) e^{-ny^2}dy
=f^{\prime}(x)+
\frac{2\sqrt{n}}{\sqrt{\pi}}\int_{0}^{n^{-1/4}} \left[f^{\prime}(y+x)-f^{\prime}(x)\right]e^{-ny^2}dy\\
+\frac{2\sqrt{n}}{\sqrt{\pi}}\int^{\infty}_{n^{-1/4}} \left[f^{\prime}(y+x)-f^{\prime}(x)\right]e^{-ny^2}dy
\end{eqnarray*}

Using the asymptotic expansion for the error function for large arguments, we have:
\begin{eqnarray*}
\left|\frac{2\sqrt{n}}{\sqrt{\pi}}\int_0^\infty f^{\prime}(y+x) e^{-ny^2}dy-f^{\prime}(x)\right|\leq
C||f^{\prime}||_{\infty}\left|n^{-1/4}+e^{-\sqrt{n}}\right|
\end{eqnarray*}

Applying this to both integrals gives us:
\begin{eqnarray*}
\frac{n^{3/2}}{\sqrt{\pi}}\int f(x)f(y)|y-x|e^{-n(x-y)^2}dydx=\frac{\sqrt{n}}{\sqrt{\pi}}\int f(x)^2 dx+O\left(n^{-1/4}\right)
\end{eqnarray*}

This is as we desired. To handle the other terms, rewrite:

\begin{eqnarray*}
\frac{4n\sqrt{n}}{\sqrt{\pi}}\int f(x)f(y)(y-x)e^{-n(x-y)^2}\\
\Big[e^{-nx^2}\frac{1}{2\sqrt{\pi}}\int_0^{ny^2}\frac{e^{-t}}{\sqrt{t}}
\sum_{m=0}^{n-1}\frac{(2\sqrt{n}x\sqrt{t})^{2m}}{(2m)!}dt]dxdy
\\
=
\frac{2\sqrt{n}}{\sqrt{\pi}}\int f(x)f(y)\frac{d}{dy}[e^{-n(x-y)^2}]\\
\Big[e^{-nx^2}\frac{\sqrt{n}}{\sqrt{\pi}}\int_0^{y}e^{-nt^2}
\sum_{m=0}^{n-1}\frac{(2nxt)^{2m}}{(2m)!}dt\Big]dxdy
\end{eqnarray*}

Integrating by parts, we may rewrite this quantity as:
\begin{eqnarray*}
-
\frac{2n}{\pi}\int f(x)f^{'}(y)e^{-n(x-y)^2}
\Big[e^{-nx^2}\int_0^{y}e^{-nt^2}
\sum_{m=0}^{n-1}\frac{(2nxt)^{2m}}{(2m)!}dt\Big]dxdy
\\
-
\frac{2n}{\pi}\int f(x)f(y)e^{-n(x-y)^2}
\Big[e^{-nx^2}e^{-ny^2}
\sum_{m=0}^{n-1}\frac{(2nxy)^{2m}}{(2m)!}\Big]dxdy
\end{eqnarray*}

We now need to handle the partial sums of the hyperbolic cosine function. We have (for real $0<|u|<1$) the asymptotic formula:
\begin{eqnarray*}
e^{-2nu}\sum_{m=0}^{2n-2}\frac{(2nu)^m}{m!}=1+O(1/n)
\end{eqnarray*}

Consequently:
\begin{eqnarray*}
\sum_{m=0}^{n-1}\frac{(2nxt)^{2m}}{(2m)!}=
\frac{1}{2}\sum_{m=0}^{2n-2}\frac{(2nxt)^{m}}{m!}+
\frac{1}{2}\sum_{m=0}^{2n-2}\frac{(-2nxt)^{m}}{m!}\\
=\frac{1}{2}
\Big(e^{2nxt}+e^{-2nxt}\Big)\Big(1+O(1/n)\Big)
\end{eqnarray*}

Plugging this in, we obtain at leading order:
\begin{eqnarray*}
-\frac{n}{\sqrt{\pi}}\int f(x)f^{\prime}(y)e^{-n(x-y)^2}
\Big[e^{-nx^2}\frac{1}{\sqrt{\pi}}\int_0^{y}e^{-nt^2}
e^{2nxt}dt\Big]dxdy
\\
-
\frac{n}{\sqrt{\pi}}\int f(x)f(y)e^{-n(x-y)^2}
\Big[e^{-nx^2}\frac{1}{\sqrt{\pi}}e^{-ny^2}
e^{2nxy}\Big]dxdy
\\
=
-\frac{n}{\pi}\int f(x)f^{\prime}(y)e^{-n(x-y)^2}
\int_0^{y}
e^{-n(x-t)^2}dtdxdy
\\
-
\frac{n}{\pi}\int f(x)f(y)e^{-2n(x-y)^2}dxdy
\end{eqnarray*}

We have that for any two test functions $f$ and $g$:
\begin{eqnarray*}
\int f(x)g(y)e^{-n(x-y)}dxdy=O(n^{-1/2})
\end{eqnarray*}

And additionally:
\begin{eqnarray*}
\sqrt{n}\int_0^y e^{-n(x-t)^2}dx\leq \sqrt{n}\int_{\mathcal{R}}e^{-nx^2}dx=\sqrt{\pi}
\end{eqnarray*}

Consequently, the first term is bounded. Moreover, by taking a Fourier transform we can write:
\begin{eqnarray*}
\frac{n}{\pi}\int f(x)f^{\prime}(y)e^{-n(x-y)^2}
\int_0^{y}
e^{-n(x-t)^2}dtdxdy=C\left(\int f(x)f^{\prime}(x)dx\right)+O\left(\frac{1}{n}\right)
\end{eqnarray*}

The first integral is zero, as can be seen by applying the fundamental theorem of calculus, $f^{\prime}f=\left(\frac{1}{2}f^2\right)^{\prime}$, and the compact support of $f$. Therefore this term is not
only bounded, it is actually $O(n^{-1})$.

The second term is also easily handled:
\begin{eqnarray*}
\frac{n}{\pi}\int f(x)f(y)e^{-2n(x-y)^2}dxdy\\
=\frac{\sqrt{n}}{2\pi}
\frac{\sqrt{n}}{\pi}\sqrt{\frac{\pi}{2n}}\int \overline{\hat{f}(s)}\hat{f}(s)e^{-s^2/8n}ds\\
=
\frac{\sqrt{n}}{\sqrt{2\pi}}\int f(x)^2(1+O(1/n))dx
\end{eqnarray*}

This concludes the proof $\P$.

Plugging these lemmas into our expression for the limiting variance, we obtain:
\begin{eqnarray*}
\mbox{Var}(S_{2n}[f])\approx
\frac{\sqrt{n}}{\sqrt{\pi}}\int f(x)^2dx-
\sqrt{n}\left(\frac{1}{\sqrt{2\pi}}-\frac{1}{\sqrt{\pi}}\right)\int f(x)^2 dx
-\frac{\sqrt{n}}{\sqrt{2\pi}}\int f(x)^2 dxdy\\
=2\left(\frac{\sqrt{n}}{\sqrt{\pi}}-\frac{\sqrt{n}}{\sqrt{2\pi}}\right)\int f(x)^2 dx
=\sqrt{n}\left(\frac{2-\sqrt{2}}{\sqrt{\pi}}\right)\int f(x)^2 dx
\end{eqnarray*}

\subsection{Higher Cumulants in the Real/Real Case: Proof of Theorem \ref{nh:line}}
In this subsection, we will show that the (normalized) higher cumulants vanish, giving a Gaussian limiting distribution and concluding the proof of Theorem \ref{nh:line}.
\begin{proposition}
Let $f$ be a  function compactly supported on the interval $(-1,1)$. Then, for $m \geq 3$, the $m$-th cumulant of the linear statistic $n^{-1/4}S_n[f]$
of the real Ginibre ensemble goes to zero as $n$ becomes large.
\end{proposition}

By the formula for cumulants of quaternion determinants, and by an argument identical to the one for higher cumulants in the complex/complex case, we know that the $N$-th cumulant of $S_n[f]$ is given by:
\begin{eqnarray*}
\mathcal{C}_{N}[f]=\int\Bigg(\sum_{V_i\in X_N}
f(x_{1})^{|V_{i,1}|}...f(x_{|V_i|})^{|V_{i,|V_i|}|}
(-1)^{|V_i|-1}\\
\times
\sum_{\sigma\in C[|V_i|]}
\left(Q_{2n}(x_{1},x_{\sigma(1)})Q_{2n}(x_{\sigma(1)},x_{\sigma^2(1)})...
Q_{2n}(x_{\sigma^{-1}(1),x_1})
\right)
\Bigg)dx_1...dx_{|V_i|}
\end{eqnarray*}

Here, we have defined the entries of the quaternion matrix $Q_n$ as:

\begin{eqnarray*}
 Q_{2n}(j,k)=S_{2n}(x_j,x_k)e_0+\frac{-i}{2}(D_{2n}(x_{j},x_{k})+I_{n}(x_j,x_k))e_1\\
+\frac{1}{2}(D_{n2}(x_{j},x_{k})-I_{n}(x_j,x_k))e_2
\end{eqnarray*}

Notice that the product $Q_{2n}(x_{1},x_{\sigma(1)})Q_{2n}(x_{\sigma(1)},x_{\sigma^2(1)})...
Q_{2n}(x_{\sigma^{-1}(1)},x_1)$ is generically a quaternion, not a scalar. However, because the determinant of a self-dual quaternion matrix is a scalar, all the quaternion terms will cancel out during summation necessarily and can therefore be ignored -- in other words, while analyizing a generic product of this form, we can throw out all terms which don't have coefficient $e_0=I$ (after the multiplication is performed, of course).

It is often easier to replace quaternion multiplication with the equivalent matrix multiplication. For instance, we can expand $Q_n(x_1,x_2)Q_n(x_2,x_1)$ using matrix multiplication as:

\[\left( \begin{array}{cc}
S_{2n}(x_1,x_2) & I_{2n}(x_1,x_2) \\
D_{2n}(x_1,x_2)& S_{2n}(x_1,x_2) \end{array}\right)
\left( \begin{array}{cc}
S_{2n}(x_2,x_1) & I_{2n}(x_2,x_1) \\
D_{2n}(x_2,x_1)& S_{2n}(x_2,x_1) \end{array}\right)\]
\[=\left( \begin{array}{cc}
S_{2n}(x_1,x_2) & I_{2n}(x_1,x_2) \\
D_{2n}(x_1,x_2)& S_{2n}(x_1,x_2) \end{array}\right)
\left( \begin{array}{cc}
S_{2n}(x_1,x_2) & -I_{2n}(x_1,x_2) \\
-D_{2n}(x_1,x_2)& S_{2n}(x_1,x_2) \end{array}\right)\]
\[= \left( \begin{array}{cc}
S_{2n}(x_1,x_2)^2-I_{2n}(x_1,x_2)D_{2n}(x_1,x_2)  & 0 \\
0 & S_{2n}(x_1,x_2)^2-I_{2n}(x_1,x_2)D_{2n}(x_1,x_2)   \end{array}\right)\]

Along with $Q_{2n}(x,x)=S_{2n}(x,x)$ gives us an alternate justification for the expression for the variance, which we had previously just read off of Pfaffians.

We would like to show that each summand corresponding to a specific $V_i\in X_N$ and specific $\sigma\in C[|V_i|$, for $N$ fixed, is small in the limit -- this will obviously imply the smallness of $\mathcal{C}_N[f]$. This program will require several steps: first, we will argue that terms which correspond to choosing the $e_0$ term from each factor during quaternion multiplication do not contribute to higher cumulants. After that, we will argue that all other terms essentially reduce to the first case after an approximation of the identity type manipulation, and the result will follow.

We begin with:

\begin{lemma}
\label{SN}
For some Schwartz function $f$ and partition $\{V_1\},...,\{V_m\}$ of $[1,....,N]$, we consider the following integral expression:
\begin{eqnarray*}
\left(\frac{(2n)^{m/2}}{n^{N/4}}\right)\int
f^{|V_1|}(x_{1})...f^{|V_m|}(x_m)
\left(S_n(\sqrt{2n}x_{1},\sqrt{2n}x_{2})...
(S_n(\sqrt{2n}x_{m},\sqrt{2n}x_{1})
\right)
dx_1...dx_{m}
\end{eqnarray*}
If $N\geq 3$, then this expression has magnitude at most $O(n^{-1/4})$.
\end{lemma}

\textbf{Proof:}

We will write $f_i(x_i)$ for $f^{|V_i|}(x_i)$. By elementary Fourier analysis:
\begin{eqnarray*}
\int f_n(x_n) e^{-n|x_{n-1}-x_{n}|^2}e^{-n|x_n-x_1|^2}dx_n
=\frac{1}{2\pi}
\int F[f_n(x_n)e^{-n|x_{n-1}-x_{n}|^2}](t_n)\overline{F[e^{-n|x_n-x_1|^2}](t_n)}dt_n\\
=\frac{1}{4\pi^2}\frac{\pi}{n}
\int \left(\int \hat{f}_n(t_n-t_{n-1})e^{-|t_{n-1}|^2/n}e^{-it_{n-1}x_{n-1}}dt_{n-1}\right)
e^{-t_n^2/4n}e^{it_nx_1}dt_n
\end{eqnarray*}

Similarly:
\begin{eqnarray*}
\int f_{n-1}(x_{n-1})e^{-n|x_{n-2}-x_{n-1}|^2}\int f_n(x_n) e^{-n|x_{n-1}-x_{n}|^2}e^{-n|x_n-x_1|^2}\\
=\frac{1}{4\pi^2}\frac{\pi}{n}\int f_{n-1}(x_{n-1})e^{-n|x_{n-2}-x_{n-1}|^2}e^{-it_{n-1}x_{n-1}}\\
\times
\left(
\hat{f}_n(t_n-t_{n-1})e^{-|t_{n-1}|^2/n}
e^{-t_n^2/4n}e^{it_nx_1}\right)dt_{n-1}dt_n
\\
=
\frac{1}{8\pi^3}\frac{\pi}{n}\sqrt{\frac{\pi}{n}}\int \hat{f}_{n-1}(t_{n-1}-t_{n-2})
e^{-|t_{n-2}|^2/4n}e^{-it_{n-2}x_{n-2}}\\
\times
\left(
\hat{f}_n(t_n-t_{n-1})e^{-|t_{n-1}|^2/n}
e^{-t_n^2/4n}e^{it_nx_1}\right)dt_{n-1}dt_n
\end{eqnarray*}

Continuing this pattern, we obtain:
\begin{eqnarray*}
(2n)^{m/2}\int f_1(x_1)f_2(x_2)...f_m(x_m)e^{-n(|x_1-x_2|^2+|x_2-x_3|^2+...+|x_{m-1}-x_m|^2+|x_m-x_1|^2)}
\prod_{k=1}^m dx_k\\
=\left(\frac{\sqrt{2n}}{2\pi}\right)^{m}\left(\frac{\pi}{n}\right)^{m/2}
\int \hat{f}_m(t_m-t_{m-1})\hat{f}_{m-1}(t_{m-1}-t_{m-2})...\hat{f}_{1}(t_{1}-t_m)\\
\times
\mbox{exp}\left(-(\sum_{j=1}^m|t_j|^2)/4n\right)\prod_{k=1}^m dt_k
\end{eqnarray*}

Using a change of variables:
\begin{eqnarray*}
(2n)^{m/2}\int f_1(x_1)f_2(x_2)...f_m(x_m)e^{-n(|x_1-x_2|^2+|x_2-x_3|^2+...+|x_{m-1}-x_m|^2+|x_m-x_1|^2)}
\prod_{k=1}^m dx_k\\
=\left(\frac{\sqrt{2n}}{2\pi}\right)^{m}\left(\frac{\pi}{n}\right)^{m/2}
\int \overline{\hat{f}_m(t_1+...+t_{m-1})}\hat{f}_{m-1}(t_{m-1})...\hat{f}_{1}(t_{1})\\
\times
\mbox{exp}\left(-|t_m|^2/4n\right)
(1+O(1/n))
\prod_{k=1}^m dt_k
\end{eqnarray*}

The integral of $\mbox{exp}\left(-|x|^2/4n\right)$ grows like $\sqrt{n}$, and by Plancharel's theorem we have that:
\begin{eqnarray*}
\int \prod_{k=1}^m f_k(x)dx=\left(\frac{1}{2\pi}\right)^m\int \overline{\hat{f}_m(t_1+...+t_{m-1})}\hat{f}_{m-1}(t_{m-1})...\hat{f}_{1}(t_{1})
\end{eqnarray*}

Consequently, the magnitude of this whole expression is at most $O(\sqrt{n})$, so if we normalize each $f_k$ by a factor of $n^{-1/4}$ and we have at least 3 such factors (which is where the assumption $m\geq 3$ comes in), the whole expression goes to zero like at most $O(n^{-1/4})$. $\P$

This lemma assures us that terms consisting soley of factors of the type $S_{2n}(x_i,x_j)$ do not contribute to the limiting behavior of the $N$-th cumulant for $N\geq 3$. Next, we extend this result to terms consisting soley of factors of the type $S_{2n}(x_i,x_j)$ 
or $D_{2n}(x_i,x_j)$ .

\begin{lemma}
\label{DN}
Consider the expression:
\begin{eqnarray*}
\frac{(2n)^{m/2}}{n^{-M/4}}\int f(x_1)...f(x_m)J^1(\sqrt{2n}x_1,\sqrt{2n}x_2)J^2(\sqrt{2n}x_2,\sqrt{2n}x_3)\\
....J^m(\sqrt{2n}x_m,\sqrt{2n}x_1)dx_1...dx_m
\end{eqnarray*}
Here, each $J^i(\cdot, \cdot)$ is either equal to $D_{2n}(\cdot, \cdot)$ or $S_{2n}(\cdot, \cdot)$. If $M\geq 3$, this entire expression is $O(n^{-1/4})$.
\end{lemma}

\textbf{Proof:}

We can expand:
\begin{eqnarray*}
\int f(x)\left|S_{2n}(\sqrt{2n}x,\sqrt{2n}y)D_{2n}(\sqrt{2n}x,\sqrt{2n}z)\right|dx
=f(y)\int  \left|S_{2n}(\sqrt{2n}x,\sqrt{2n}y)D_{2n}(\sqrt{2n}x,\sqrt{2n}z)\right| dx\\
+\int [f(x)-f(y)]\times \left|S_{2n}(\sqrt{2n}x,\sqrt{2n}y)D_{2n}(\sqrt{2n}x,\sqrt{2n}z)\right|dx
\end{eqnarray*}

The second term is small, as we can use $|S_{2n}|\leq O(1)$ to estimate:
\begin{eqnarray*}
\int_{|x-y|\geq n^{-1/4}} [f(x)-f(y)] \left|S_{2n}(\sqrt{2n}x,\sqrt{2n}y)D_{2n}(\sqrt{2n}x,\sqrt{2n}z)\right| dx
\leq C e^{-\sqrt{n}}
\end{eqnarray*}

And:
\begin{eqnarray*}
\int_{|x-y|\leq n^{-1/4}} [f(x)-f(y)] \left|S_{2n}(\sqrt{2n}x,\sqrt{2n}y)D_{2n}(\sqrt{2n}x,\sqrt{2n}z)\right|dx
\leq C ||f^{\prime}||_{\infty}n^{-1/4}
\end{eqnarray*}

It remain to control the first term:
\begin{eqnarray*}
\sqrt{n}\int \left|S_{2n}(\sqrt{2n}x,\sqrt{2n}y)D_{2n}(\sqrt{2n}x,\sqrt{n}z)\right|dx
=n\frac{1}{2\pi} \int e^{-n(x-y)^2}e^{-n(x-z)^2}|x-z|dx\\
=n\frac{1}{2\pi}  e^{-n(y-z)^2/2}\int |x|e^{-n(\sqrt{2}x-(y-z)/\sqrt{2}}dx
\end{eqnarray*}

Integrating, this last display is just:
\begin{eqnarray*}
-\frac{2}{8\pi} e^{-n(y-z)^2/2}\left(\sqrt{\pi}\frac{(y-z)}{\sqrt{2}}\mbox{erf}\left(\frac{(y-z)}{\sqrt{2}}-x\right)
+\frac{1}{2\pi} e^{-n(x-(y-z)/\sqrt{2})^2} \right)\Bigg|_{0}^\infty\\
=\frac{1}{4\pi} -e^{-n(y-z)^2/2}\sqrt{\pi}\frac{(y-z)}{\sqrt{2}}\left[\sqrt{\pi}-\mbox{erf}\left(\frac{(y-z)}{\sqrt{2}}\right)\right]
+\frac{1}{4\pi} e^{-n(y-z)^2}
\end{eqnarray*}
 
We have the estimate $\left|u\exp(-nu^2)\right|\leq O(n^{-1/4})$, by considering small and large $u$ seperately, so we can write (up to small correction terms which vanish in the limit):
\begin{eqnarray*}
\sqrt{n}\int \left|f(x)S_{2n}(\sqrt{2n}x,\sqrt{2n}y)D_{2n}(\sqrt{2n}x,\sqrt{n}z)\right| dx\leq \frac{1}{\sqrt{4\pi}}\left|f(z) S_{2n}(\sqrt{2n}y,\sqrt{2n}z)\right|
\end{eqnarray*}

And by a similar argument:
\begin{eqnarray*}
\sqrt{n}\int\left|f(x) D_{2n}(\sqrt{2n}x,\sqrt{2n}y)D_{2n}(\sqrt{2n}x,\sqrt{n}z)\right| dx\leq \frac{1}{\sqrt{4\pi}}\left|f(z)D_{2n}(\sqrt{2n}y,\sqrt{2n}z)\right|
\end{eqnarray*}

The point of these estimate  is that we may now eliminate terms of the form $D_{2n}$. For instance, we have:
\begin{eqnarray*}
\int f(x)f(y)D_{2n}(\sqrt{2n}x,\sqrt{2n}y)D_{2n}(\sqrt{2n}x,\sqrt{2n}y)dxdy
=\int f(x)f(y)\left|D_{2n}(\sqrt{2n}x,\sqrt{2n}y)\right|^2dxdy\\
=\frac{2n}{2\pi}\int f(x)\left(\int f(y)|x-y|^2e^{-2n(x-y)^2} dy\right) dx
\approx
\frac{n}{\pi}\frac{\sqrt{\pi}}{4n\sqrt{2n}}\int f(x)^2 dx
\end{eqnarray*}

We may conclude that for $M\geq 3$:
\begin{eqnarray*}
\frac{(\sqrt{2n})^2}{n^{-M/4}}\int f(x)f(y)D_{2n}(\sqrt{2n}x,\sqrt{2n}y)D_{2n}(\sqrt{2n}x,\sqrt{2n}y)dxdy
\leq O(n^{-1/4})
\end{eqnarray*}

Consequently, using the identities we have derived:
\begin{eqnarray*}
\left|\frac{(2n)^{m/2}}{n^{-M/4}}\int f(x_1)....f(x_m)D_{2n}(\sqrt{2n}x_1,x_2)....D_{2n}(\sqrt{2n}x_m,\sqrt{2n}x_1)dx_1...dx_m\right|\\
\leq C\frac{n}{n^{-M/4}} \int \left|f(x)f(y)\right|\times\left|D_{2n}(\sqrt{2n}x,\sqrt{2n}y)\right|^2dxdy \leq O(n^{-1/4})
\end{eqnarray*}

If instead of some of the $D_{2n}(\sqrt{2n}x,\sqrt{2n}y)$ terms were  $S_{2n}(\sqrt{2n}x,\sqrt{2n}y)$, we could still use the same trick to eliminate all of the $D_{2n}$ terms, 
and then apply Lemma \ref{SN}. For instance:
\begin{eqnarray*}
\frac{(2n)^{3/2}}{n^{-M/4}}\left|\int f(x)f(y)f(z)  S_{2n}(\sqrt{2n}x,\sqrt{2n}y) S_{2n}(\sqrt{2n}y,\sqrt{2n}z) D_{2n}(\sqrt{2n}z,\sqrt{2n}x)dz dy dx\right|\\
\leq C\frac{2n}{n^{-M/4}}\int \left|f(x)^2f(y)\right|\times \left| S_{2n}(\sqrt{2n}x,\sqrt{2n}y) S_{2n}(\sqrt{2n}y,\sqrt{2n}x)\right|dy dx
\end{eqnarray*}

By Lemma \ref{SN}, this goes to zero. Generic terms are handled the exact same way. $\P$

It remains to consider expressions involving terms of the form $I_{2n}(x_i,x_j)$. To do this, we will need to investigate the structure of $\mathcal{C}_{N}[f]$ more closely. We have the following easy lemma:

\begin{lemma}
Write the product $Q_{2n}(x_{1},x_{\sigma(1)})Q_{2n}(x_{\sigma(1)},x_{\sigma^2(1)})...
Q_{2n}(x_{\sigma^{-1}(1),x_1})$ in the form of a matrix. Then each summand in the top left and bottom right entry has an equal number of terms of the form $I_n(x_i,x_j)$ as terms of the form $D_n(x_i,x_j)$ (this number may change from summand to summand). 

Moreover, they are adjacent -- we can take each $I_n(x_i,x_j)$ to be preceded by a $D_n(x_k,x_i)$.
\end{lemma}

\textbf{Proof:}

Consider a matrix $Q$ whose entries are elements $S$, $I$ and $D$ which reside in some non-commutative algebra.
\[Q= \left( \begin{array}{cc}
S & I \\
D & S \end{array}\right)\]

We will investigate the product of $n$ such matrices:
\[\prod_{i=1}^n Q= \left( \begin{array}{cc}
S & I \\
D & S \end{array}\right)...\left( \begin{array}{cc}
S & I \\
D & S \end{array}\right)\]

We will begin by proving two properties of $Q^n$, using induction. First, for $n=1$, it is obvious that the two diagonal entries and two off-diagonal entries are identical except for a swapping of $I$'s for $D$'s. Assume then that this property holds for the product of $Q$ with itself $n-1$ times, $Q^{n-1}$, and expand $Q^{n}=Q^{n-1}Q$ as:
\[= \left( \begin{array}{cc}
W & X \\
Y & Z \end{array}\right)\left( \begin{array}{cc}
S & I \\
D & S \end{array}\right)\]
\[\left( \begin{array}{cc}
WS+XD & WI+XS \\
YS+ZD & YI+ZS \end{array}\right)\]

If we swap $I$ with $D$ and use the induction assumption (which is that this swap also interchanges $W$ for $Z$ and $X$ and $Y$), we obtain:
 \[= \left( \begin{array}{cc}
ZS+YI & DZ+YS \\
XS+WI & XD+WS \end{array}\right)\]

Therefore $Q^n$ has this property also, and so by induction $Q^m$ has this property for all $m$.

Consider now the effect of matrix multiplication on the top left and top right entries of $Q$. This effect is captured in the following mapping:
\begin{eqnarray*}
(X,Y)\to (XS+YD,XI+YS)
\end{eqnarray*}

Here, the ordered pair $(\cdot,\cdot)$ represents the top left and top right entries of a matrix.

For $n=1$, it is clear that the top left entry of $Q$ has the same number of $I$'s and $D$'s (in this case, zero) and that the top right entry has one more $I$ than $D$ (in this case, one to zero). Assume that this property holds for $Q^{n-1}$. The, since we are assuming that $X$ has the same number of $I$'s and $D$'s, $XS$ does also. Similarly, since $Y$ has one extra $I$, $YD$ has the same number of each again, and therefore $XS+YD$ does also. It follows just as easily that $XI$ and $YS$ both have one more $I$ than $D$ (just using the assumptions again), and therefore $Q^{n+1}$ has the same number of $I$'s and $D$'s in its top left entry and one more $I$ than $D$ in its top right entry. By induction, all products $Q^m$ have this property.

Combining both observations and identifying $Q^m$ with $\prod_{i=1}^m Q(x_{i},x_{i+1})$ (with $x_{m+1}=x_{1}$) in the obvious manner, the lemma is proven.   $\P$

This lemma is exactly what we need to prove:

\begin{lemma}
Consider the expression:
\begin{eqnarray*}
\frac{(2n)^{m/2}}{n^{-M/4}}\int f(x_1)...f(x_m)J^1(\sqrt{2n}x_1,\sqrt{2n}x_2)J^2(\sqrt{2n}x_2,\sqrt{2n}x_3)\\
....J^m(\sqrt{2n}x_m,\sqrt{2n}x_1)dx_1...dx_m
\end{eqnarray*}
Here, each $J^i(\cdot, \cdot)$ is either equal to $D_{2n}(\cdot, \cdot)$ , $S_{2n}(\cdot, \cdot)$, or $I_{2n}(\cdot, \cdot)$ (subject to the condition that all $I_{2n}$ are immediately preceded by a $D_{2n}$) . If $M\geq 3$, this entire expression is $O(n^{-1/4})$.
\end{lemma}

\textbf{Proof:}

If none of the $J^i$ are equal to $I_{2n}$, then this follows from our previous results. If on the other hand some $J^{i}$ is equal to $I_{2n}$, then the previous lemma assures us that $J^{i-1}$ is equal to $D_{2n}$. 
We therefore need to investigate terms of the form $D_{2n}(\cdot, \cdot)I_{2n}(\cdot, \cdot)$. 

Let us define:
\begin{eqnarray*}
\tilde{I}_{2n}(x,y)=
e^{-x^2/2}\frac{1}{2\sqrt{\pi}}\int_0^{y^2/2}\frac{e^{-t}}{\sqrt{t}}
\sum_{m=0}^{n-1}\frac{(\sqrt{2t}x)^{2m}}{(2m)!}dt\\-
e^{-y^2/2}\frac{1}{2\sqrt{\pi}}\int_0^{x^2/2}\frac{e^{-t}}{\sqrt{t}}
\sum_{m=0}^{n-1}\frac{(y\sqrt{2t})^{2m}}{(2m)!} dt
\end{eqnarray*}

Then we have:
\begin{eqnarray*}
I_{2n}(x,y)=\tilde{I}_{2n}(x,y)+\frac{1}{2}\mbox{sgn}(x-y)
\end{eqnarray*}

Now, recall that we have the identity:

\begin{eqnarray*}
e^{-nx^2}\frac{1}{2\sqrt{\pi}}\int_0^{ny^2}\frac{e^{-t}}{\sqrt{t}}
\sum_{m=0}^{n-1}\frac{(\sqrt{2t}\sqrt{2n}x)^{2m}}{(2m)!}dt\\
=e^{-nx^2}\frac{\sqrt{n}}{\sqrt{\pi}}\int_0^{y}e^{-nt^2}
\sum_{m=0}^{n-1}\frac{(2ntx)^{2m}}{(2m)!}dt\
\end{eqnarray*}

And:
\begin{eqnarray*}
\sum_{m=0}^{n-1}\frac{(2nxt)^{2m}}{(2m)!}
=\frac{1}{2}
\Big(e^{2nxt}+e^{-2nxt}\Big)\Big(1+O(1/n)\Big)
\end{eqnarray*}

Combining, we have:
\begin{eqnarray*}
\tilde{I}_{2n}(x,y)\approx\frac{\sqrt{n}}{2\sqrt{\pi}}\int_0^{y}e^{-n(x-t)^2}dt-
\frac{\sqrt{n}}{2\sqrt{\pi}}\int_0^{x}e^{-n(y-t)^2}dt
\end{eqnarray*}

We will also need the derivative of the right hand side:
\begin{eqnarray*}
\frac{d}{dx}\left(\frac{\sqrt{n}}{2\sqrt{\pi}}\int_0^{y}e^{-n(x-t)^2}dt-
\frac{\sqrt{n}}{2\sqrt{\pi}}\int_0^{x}e^{-n(y-t)^2}dt \right)\\
=\frac{\sqrt{n}}{2\sqrt{\pi}}\left(\int_0^{y}-\frac{d}{dt}\left(e^{-n(x-t)^2}\right)dt-e^{-n(y-x)^2}\right)\\
=\frac{\sqrt{n}}{2\sqrt{\pi}}e^{-nx^2}-\frac{\sqrt{n}}{\sqrt{\pi}}e^{-n(y-x)^2}
\end{eqnarray*}

Then we can write, using integration by parts:
\begin{eqnarray*}
2n\int f(x)(x-z)e^{-n(z-x)^2}\tilde{I}_{2n}(\sqrt{2n}x,\sqrt{2n}y)dx\\
=-\int f^{\prime}(x)e^{-n(z-x)^2}\tilde{I}_{2n}(\sqrt{2n}x,\sqrt{2n}y)dx\\
-\int f(x)e^{-n(z-x)^2}\left(\frac{\sqrt{n}}{2\sqrt{\pi}}e^{-nx^2}-\frac{\sqrt{n}}{\sqrt{\pi}}e^{-n(y-x)^2} \right)dx
\end{eqnarray*}

By a Fourier transform:
\begin{eqnarray*}
\frac{\sqrt{n}}{\sqrt{\pi}}\int f(x)e^{-n(y-x)^2}dx=\frac{1}{2\pi}\frac{\sqrt{n}}{\sqrt{\pi}}\sqrt{\frac{\pi}{n}}\int \hat{f}(\xi)e^{-\xi^2/4n}e^{i\xi y}d\xi\\
=\frac{1}{2\pi}\int \hat{f}(\xi)e^{i\xi y}\left(1+O(1/n)\right) d\xi\\
=f(y)(1+O(1/n))
\end{eqnarray*}

And consequently we have:
\begin{eqnarray*}
\sqrt{2n}\int f(x)D_{2n}(\sqrt{2n}z,\sqrt{2n}x)\tilde{I}_{2n}(\sqrt{2n}x,\sqrt{2n}y)dx\\
=-\frac{1}{2\sqrt{2\pi}}f(z)e^{-nz^2}+\frac{1}{\sqrt{2\pi}}\frac{1}{\sqrt{2}}e^{-n(z-y)^2/2}f\left(\frac{y+z}{2}\right)+O\left(\frac{1}{\sqrt{n}}\right)
\end{eqnarray*}

The first term on the right hand side is small (this can be seen either by taking $|x|\geq n^{-1/4}$ without loss of generality, as all of our integrals are bounded, or by appealing to the bound attained in Lemma \ref{SN} directly). Consequently, we have:
\begin{eqnarray*}
\sqrt{2n}\int f(x)D_{2n}(\sqrt{2n}z,\sqrt{2n}x)\tilde{I}_{2n}(\sqrt{2n}x,\sqrt{2n}y)dx\approx 
\frac{1}{2\sqrt{\pi}}e^{-n(z-y)^2/2}f\left(\frac{y+z}{2}\right)
\end{eqnarray*}

Next, we consider:
\begin{eqnarray*}
\sqrt{2n}\int f(x)D_{2n}(\sqrt{2n}z,\sqrt{2n}x)\frac{1}{2}\mbox{sgn}(x-y)dx=
2n\int f(x)\frac{1}{\sqrt{2\pi}}(x-z)e^{-n(z-x)^2}\frac{1}{2}\mbox{sgn}(x-y)dx\\
=\frac{n}{\sqrt{2\pi}}\int_y^\infty f(x)(x-z)e^{-n(z-x)^2}dx
-\frac{n}{\sqrt{2\pi}}\int_{-\infty}^y f(x)(x-z)e^{-n(z-x)^2}dx\\
=\frac{1}{2\sqrt{2\pi}}\left(\int_y^\infty f(x)\frac{d}{dx}\left[e^{-n(z-x)^2}\right]dx
-\int_{-\infty}^y f(x)\frac{d}{dx}\left[e^{-n(z-x)^2}\right]dx\right)
\end{eqnarray*}

Integrating by parts, this is:
\begin{eqnarray*}
\frac{-1}{2\sqrt{2\pi}}\int_y^\infty f^{\prime}(x)\mbox{sgn}(x-y)e^{-n(z-x)^2}dx
-\frac{1}{\sqrt{2\pi}}f(y)e^{-n(z-y)^2}
\end{eqnarray*}

And we conclude:
\begin{eqnarray*}
\sqrt{2n}\int f(x)D_{2n}(\sqrt{2n}z,\sqrt{2n}x)\frac{1}{2}\mbox{sgn}(x-y)dx\\
=\frac{-1}{\sqrt{2\pi}}f(y)e^{-n(z-y)^2}+O\left(\frac{1}{\sqrt{n}}\right)
\end{eqnarray*}

We have therefore shown that, since $||f||_{\infty}<\infty$, up to small terms which vanish in the limit the following estimate holds:
\begin{eqnarray*}
\sqrt{2n}\left|\int f(x)D_{2n}(\sqrt{2n}z,\sqrt{2n}x)I_{2n}(\sqrt{2n}x,\sqrt{2n}y)dx\right|\\
\leq C\left(e^{-n(z-y)^2/2}+e^{-n(z-y)^2}\right)
\end{eqnarray*}

Applying this argument to all $I_{2n}$ terms which appear, and then applying the argument from the proof of Lemma \label{DN} to handle all of the remaining $D_{2n}$ terms (if any), we see that
it suffices to prove the lemma when each $J^{i}(\sqrt{2n}x,\sqrt{2n}y)$ is either of the form $e^{-n(x-y)^2/2}$ or $S_{2n}(\sqrt{2n}x,\sqrt{2n}y)$. However, because $S_{2n}(\sqrt{2n}x,\sqrt{2n}y)$ is proportional to
$e^{-n(x-y)^2}$, this shouldn't prove too difficult.

Indeed, in the proof of Lemma \ref{SN}, we demonstrated that as long as $M \geq 3$:
\begin{eqnarray*}
\frac{(2n)^{m/2}}{n^{M/4}}\int f_1(x_1)f_2(x_2)...f_m(x_m)e^{-n(|x_1-x_2|^2+...+|x_m-x_1|^2)}
\prod_{k=1}^m dx_k
\leq O\left(\frac{1}{n^{1/4}}\right)
\end{eqnarray*}

Inspecting the proof, we see that nothing except the implied constant would change if we swapped one or more of the $|x_i-x_j|^2$ terms in the exponent with $|x_i-x_j|^2/2$, and therefore the lemma is proven $\P$.

\end{document}